\documentclass[leqno,12pt]{article}

\newlength{\textlarg}

\usepackage{amssymb, ulem}
\usepackage{euscript}
\usepackage{mathrsfs} 
\usepackage[latin1]{inputenc}
\usepackage{amsmath}
\usepackage{amsthm}
\usepackage{amsfonts}
\usepackage{enumerate}

\usepackage[colorlinks = true, linkcolor = blue, urlcolor  = blue, citecolor = blue, anchorcolor = blue]{hyperref} 

 \usepackage{graphicx}
 \usepackage{float} 

\evensidemargin\oddsidemargin

\def\thesection{\arabic{section}}
\renewcommand{\theequation}{\thesection.\arabic{equation}}

\newtheorem{theorem}{Theorem}[section]
\newtheorem{lemma}[theorem]{Lemma}
\newtheorem{proposition}[theorem]{Proposition}
\newtheorem{corollary}[theorem]{Corollary}

\newtheorem{notation}[theorem]{Notation}

\theoremstyle{definition}   
\newtheorem{remark}[theorem]{Remark}

\newcommand{\eqnsection}{
\renewcommand{\theequation}{\thesection.\arabic{equation}}
    \makeatletter
    \csname  @addtoreset\endcsname{equation}{section}
    \makeatother}
\eqnsection

\def\r{{\mathbb R}}
\def\e{{\mathbb E}}
\def\p{{\mathbb P}}

\def\ee{\mathrm{e}}
\def\d{\mathrm{d}}
\def\law{\, {\buildrel \scriptsize{\rm (law)} \over =}\, }

\begin{document}


 \vglue30pt

\centerline{\large\bf   Path decompositions of perturbed reflecting Brownian motions}

\bigskip
\bigskip

 \centerline{by}

\medskip

 \centerline{Elie A\"{i}d\'ekon\footnote{\scriptsize Fudan University, LPSM, Sorbonne Universit\'e Paris VI, and Institut Universitaire de France, {\tt elie.aidekon@upmc.fr}},   Yueyun Hu\footnote{\scriptsize LAGA, Universit\'e Sorbonne Paris Nord, France, {\tt yueyun@math.univ-paris13.fr}}, 
and Zhan Shi\footnote{\scriptsize LPSM, Sorbonne Universit\'e Paris VI, France, {\tt zhan.shi@upmc.fr}}}

\bigskip
\bigskip

\centerline{\it In honour of Professor Ron Doney}

\centerline{\it on the occasion of his 80th birthday}

\bigskip

{\leftskip=2truecm \rightskip=2truecm \baselineskip=15pt \small

\noindent{\slshape\bfseries Summary.} We are interested in path decompositions of a perturbed reflecting Brownian motion (PRBM) at the hitting times and at the minimum. Our study relies on the loop soups developed by Lawler and Werner~\cite{lawler-werner} and  Le Jan~\cite{lejan}-\cite{lejan_saintflour}, in particular on a result discovered by Lupu~\cite{lupu} identifying the law of the excursions of the PRBM above its past minimum with the loop measure of Brownian bridges.

\medskip

\noindent{\slshape\bfseries Keywords.} Perturbed reflecting Brownian motion, path decomposition, Brownian loop soup, Poisson--Dirichlet distribution. 

\medskip
 
\noindent{\slshape\bfseries 2010 Mathematics Subject
Classification.} 60J65.

} 

\bigskip
\bigskip

\section{Introduction}

Let $(B_t, \, t\ge 0)$ be a standard one-dimensional Brownian motion. Let
$$
{\mathfrak L}_t:= \lim_{\varepsilon\to 0} \frac{1}{\varepsilon} \int_0^t  1_{\{ 0<B_s \le \varepsilon\} } \, \mathrm{d} s \, ,
\qquad
\hbox{\rm a.s.},
$$

\noindent be the local time at time $t$ and position $0$. We take a continuous version of $({\mathfrak L}_t, \, t\ge 0)$. Let $\mu \in \r \backslash \{ 0\}$ be a fixed parameter. Consider the perturbed reflecting Brownian motion (PRBM)
\begin{equation} \label{def-muprocess}
X_t := |B_t| - \mu {\mathfrak L}_t, \qquad t\ge 0.
\end{equation}

\noindent The PRBM family contains two important special members: Brownian motion ($\mu=1$; this is seen using L\'evy's identity), and the three-dimensional Bessel process ($\mu=-1$; seen by means of L\'evy's and Pitman's identities). 

The PRBM, sometimes also referred to as the $\mu$-process and appearing in the literature as the limiting process in the winding problem for three-dimensional Brownian motion around lines (Le Gall and Yor~\cite{legall-yor87}), turns out to have remarkable properties such as the Ray--Knight theorems (Le Gall and Yor~\cite{legall-yor}, Werner~\cite{werner95}, Perman~\cite{perman}, Perman and Werner~\cite{perman-werner}), and L\'evy's arc sine law (Petit~\cite{petit}, Carmona, Petit and Yor~\cite{carmona-petit-yor94b}). The process can also be viewed as (non-reflecting) Brownian motion perturbed by its one-sided maximum (Davis~\cite{davis99}, Perman and Werner~\cite{perman-werner}, Chaumont and Doney~\cite{chaumont-doney}). As explained on page 100 of Yor~\cite{yor92}, these simply formulated and beautiful results were proved for the PRBM because of the scaling property and the strong Markov property of $(|B|, \, {\mathfrak L})$ via excursion theory.

We study in this paper path decompositions of the PRBM. Apart from their own interests, these decompositions  can be used to understand  the dual of general Jacobi stochastic flows which will be given in a forthcoming work, extending the work of Bertoin and Le Gall~\cite{blg2} who proved that the Jacobi flows of parameters $(0, \, 0)$ and $(2, \, 2)$ are dual with each other. These stochastic flows are connected to other important probabilistic objets in the study of population genetics such as flows of Fleming--Viot processes. Technically, our study of the PRBM often relies on the powerful tool of loop soups (Lawler and Werner~\cite{lawler-werner}, Le Jan~\cite{lejan}-\cite{lejan_saintflour}), and in particular, on a result discovered by Lupu~\cite{lupu} identifying the law of the excursions of the PRBM above its past minimum with the loop measure of Brownian bridges. The point of view via loop soups has two main advantages: (i) its nice properties under rerooting allows  to  shed light on or extend   some previously known results,  see Proposition \ref{p:Jjump} or Theorem \ref{p:scaling1'}, 
and  (ii) thanks to the independence structure in the Poisson point process representation of the loop measure, it helps to make arguments of conditioning rigorous:  for instance,  we show in Lemma \ref{l:perman} a  path decomposition   for PRBM, originated from Perman \cite{perman}.  

Our path decompositions focus on two families of random times of a recurrent PRBM: first hitting times (Section \ref{s:decomposition}), and times at which the PRBM reaches its past minimum (Section \ref{s:decompositionminimum}). To illustrate the kind of results we have obtained, let us state two examples. The first one, Theorem \ref{p:hitting}, yields in the special case $\mu=1$ the classical Williams' Brownian path decomposition theorem (Revuz and Yor~\cite{revuz-yor}, Theorem VII.4.9). The second, Theorem \ref{p:decomposition}, in the special case $\mu=1$, states as follows:

\medskip

{\bf Theorem \ref{p:decomposition} (special case $\mu=1$).} {\it Consider $B$ up to the first time   that ${\mathfrak L}$ reaches $1$, time-changed to remove its excursions above zero. We decompose this path at the minimum into the post- and time reversed pre- minimum processes. The two processes are independent and distributed as $X^1$ and $X^2$,  time-changed to remove their excursions above level $H^{1, 2}$, where $X^1$ and $X^2$ are independent three-dimensional Bessel processes and $H^{1, 2}$ is the last level at which the sum of the total local times of $X^1$ and $X^2$ equals $1$.} 

\medskip

The rest of this paper is organized as follows. 

$\bullet$    Section \ref{s:pre}:   we  recall Lupu \cite{lupu}'s connection between the PRBM and the Brownian loop soup, and   some known results on the PRBM. We also  study the  minimums of the PRBM considered up to its inverse  local times   (the process  $J$ defined in \eqref{def-Jx}). The main (new) result in this section is Proposition \ref{p:Jjump}, a description of the jumping times of $J$;

\begin{figure}[h]
\centering
\includegraphics[width=0.8\textwidth]{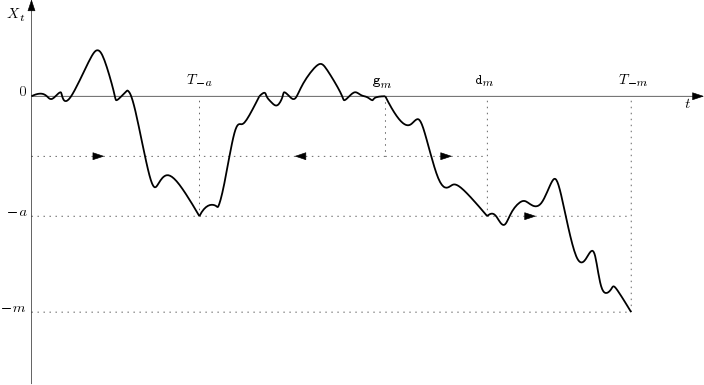}
\caption{\leftskip=1.8truecm \rightskip=1.8truecm  The four processes in Theorem \ref{p:hitting}.}

 \label{f:fig1}
\end{figure}

$\bullet$  Section \ref{s:decomposition}:  we study the path decomposition at the hitting time  of the PRBM. We prove that conditioned on its minimum, the PRBM can be split into four independent processes  (see Figure \ref{f:fig1}) and describe their laws in Theorems \ref{p:hitting} and \ref{p:hitting2};

\begin{figure}[htb]
\centering
\includegraphics[width=0.8\textwidth]{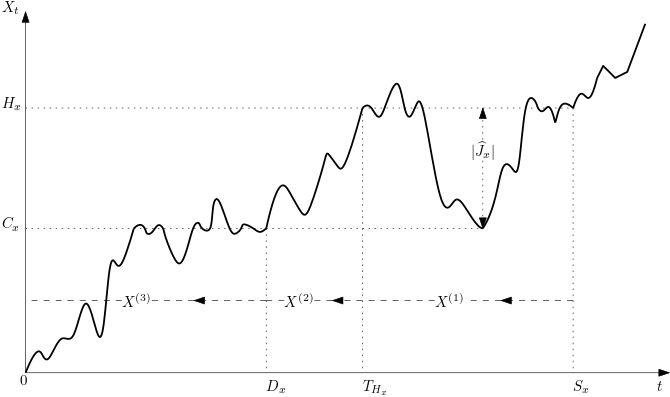}
\caption{\leftskip=1.8truecm \rightskip=1.8truecm   Under $\p^{(-\delta)}$, $X_t \to \infty$ a.s.   }

 \label{f:fig2}
\end{figure}

$\bullet$  Section \ref{s:decompositionminimum}:  we study the path decomposition   at the   minimum  of the PRBM considered up to its inverse  local time.   Theorem \ref{p:decomposition} deals with the recurrent case and  describes the laws of  the post- and time reversed pre- minimum processes.    A similar decomposition is obtained in Proposition \ref{p:mucond}  for the  transient case (see Figure \ref{f:fig2});

$\bullet$  Section \ref{s:bessel}: we extend  the perturbed Bessel process studied in Doney, Warren and Yor \cite{dwy98}  to the perturbed Bessel process with a positive local time at $0$. The main result in this section (Theorem \ref{p:scaling1'})  gives  an extension of   Theorem 2.2 of  Doney, Warren and Yor \cite{dwy98}.

\section{Preliminaries}\label{s:pre}

This section is divided into three subsections.  We recall in  Section \ref{s:loopsoup}  Lupu \cite{lupu}'s description (Proposition \ref{p:loopsoup}) on the excursions of the PRBM above its past minimum in terms  of Brownian loop soup,  and  we collect the intensities of various Poisson point processes in Lemma \ref{l:intensity}. In Section \ref{s:poisson}, we study the minimum process $J(x)$ defined in \eqref{def-Jx} and  describe the jump times of $J(\cdot)$ by means  of Poisson--Dirichlet distributions in Proposition \ref{p:Jjump}.  Finally in Section \ref{s:known} we recall some known results on the PRBM.   We also introduce some notations (in particular Notations \ref{n:pdelta} and \ref{n:<h}) which are used throughout the paper.

\subsection{The Brownian loop soup}\label{s:loopsoup}

  Lupu \cite{lupu} showed a connection between perturbed reflecting Brownian motions and the Brownian loop soup. We rely on \cite{lupu} and review this connection in this subsection.   Let ${\mathcal K}$ denote the set of continuous functions $\gamma: [0, T(\gamma)] \to \r$ with some $T(\gamma) \in (0, \infty)$, endowed with a metric $d_{\mathcal K}(\gamma, \widehat  \gamma):= |\log T(\gamma)-\log T(\widehat  \gamma)| + \sup_{0\le s\le 1} | \gamma(s T(\gamma))- \widehat \gamma(s T(\widehat \gamma))|$ for any $\gamma, \widehat \gamma \in {\mathcal K}$.  A rooted loop  is  an element $\gamma$ of ${\mathcal K}$  such that $\gamma(0)=\gamma(T(\gamma))$ (Section 3.1, p. 29). On the space of rooted loops, one defines  the measure (Definition 3.8, p.37)
$$
\mu_{\rm loop}({\rm d}\gamma) := \int_{t>0} \int_{x\in \mathbb{R}} P^t_{x,x}({\rm d} \gamma) p_t(x,x) {\rm d } x {{\rm d} t\over t},
$$
where $P^t_{x,x}$ is the distribution of the Brownian bridge of length $t$ from $x$ to $x$, and $p_t(x,x)$ is the heat kernel  $p_t(x,x) = {1\over \sqrt{2\pi t}}$. An unrooted loop is the equivalence class of all loops obtained from one another by time-shift, and $\mu_{\rm loop}^*$ denotes the projection of $\mu_{\rm loop}$ on the space of unrooted loops.  For  any fixed $\beta > 0$, the Brownian loop soup of intensity measure $\beta$ is the Poisson point process on the space of unrooted loops with intensity measure given by $  \beta \mu_{\rm loop}^*$ (Definition 4.2, p. 60). We denote it by ${\mathcal L}_\beta$.

For any real $q$, we let $\gamma-q$ denote the loop $(\gamma(t)-q,\, 0\le t \le T(\gamma))$. We write $\min \gamma$, resp. $\max \gamma$ for the minimum, resp. maximum of $\gamma$.  If $\gamma$ denotes a loop, the loop $\gamma$ rooted at its minimum is the rooted loop obtained by shifting the starting time of the loop to the hitting time of $\min \gamma$. Similarly for the loop $\gamma$ rooted at its maximum. By an abuse of notation, we will often write $\gamma$ for its range. For example $0\in \gamma$ means that $\gamma$ visits the  point $0$.

Similarly to Lupu \cite{lupu}, Section 5.2, define
\[
{\mathcal Q}^\uparrow_\beta  := \{ \min \gamma,\, \gamma \in {\mathcal L}_\beta \}, \qquad 
{\mathcal Q}^\downarrow_\beta  := \{ \max \gamma,\, \gamma \in {\mathcal L}_\beta \}.
\]

\noindent For any $q\in {\mathcal Q}^\uparrow_\beta$ and $\gamma \in {\mathcal L}_\beta$ such that $\min \gamma=q$, define ${\mathfrak e}_q^\uparrow$ as the loop $\gamma-q$  rooted at its minimum. It is an excursion above $0$. Define similarly, for any $q\in {\mathcal Q}^\downarrow_\beta$, ${\mathfrak e}_q^\downarrow$ as the excursion below $0$ given by  $\gamma-q$  rooted at  its maximum.  The point measure $\{(q,{\mathfrak e}_q^\downarrow),\, q \in {\mathcal Q}^{\downarrow}_\beta\}$ has the same distribution  as  $\{ - (q,{\mathfrak e}_q^\uparrow),\, q \in {\mathcal Q}^{\uparrow}_\beta\}$.

 \begin{notation} \label{n:pdelta}
 For $\delta>0$, let  $\p^\delta$  (resp. $\p^{(-\delta)}$) be the probability measure under which $(X_t)_{t\ge0}$ is distributed as the  PRBM  $(|B_t|- \mu {\mathfrak L}_t)_{t\ge 0}$ defined in \eqref{def-muprocess} with $\mu= \frac2\delta$ (resp. $\mu=-\frac2\delta$).   
  \end{notation}

Note  that under $\p^\delta$, $(X_t)_{t\ge0}$ is recurrent whereas under $\p^{(-\delta)}$,  $\lim_{t\to\infty} X_t =+\infty$ a.s.  Define 
 $$
 I_t := \inf_{0\le s\le t} X_s,\, t\ge 0.
 $$
 
\noindent We will use the same notation $\{ (q, {\mathfrak e}_{X,q}^\uparrow),\, q \in {\mathcal Q}^\uparrow_X\}$  to denote
\begin{itemize}
\item under $\p^\delta$: the excursions   away from  $\r \times \{0\} $  ($q$ is seen as a real number) of the process $ (I_t,X_t-I_t)$;
\item under $\p^{(-\delta)}$:  the excursions   away from  $\r \times \{0\} $ of the process $ (\widehat I_t,X_t-\widehat I_t)$ where $\widehat I_t:=\inf_{s\ge t} X_s$ . 
\end{itemize}

 The following proposition is for example Proposition 5.2 of \cite{lupu}. One can also see it from \cite{legall-yor} or \cite{AM} (together with Proposition 3.18 of \cite{lupu}).

\begin{proposition} [Lupu  \cite{lupu}] \label{p:loopsoup}
Let $\delta>0$. The point measure $\{ (q, {\mathfrak e}_{X,q}^\uparrow),\, q \in {\mathcal Q}^\uparrow_X\}$ is distributed under $\p^\delta$, respectively $\p^{(-\delta)}$, as $\{ (q, {\mathfrak e}_q^\uparrow),\, q \in {\mathcal Q}^\uparrow_{\delta \over 2} \cap (-\infty,0)\}$, resp. $\{ (q, {\mathfrak e}_{q}^\uparrow),\, q \in {\mathcal Q}^\uparrow_{\delta \over 2} \cap (0,\infty)\}$. \end{proposition}

Actually, Proposition 5.2 \cite{lupu} states the previous proposition in a slightly different way. In the same way that standard Brownian motion can be constructed from its excursions away from $0$, Lupu shows that one can construct the perturbed Brownian motions from the Brownian loop soup by ``gluing'' the loops of the Brownian loop soup rooted at their minimum and  ordered by decreasing minima.

We close this  section by collecting the intensities of various Poisson point processes. It comes from computations of \cite{lupu}.

Denote  by ${\mathfrak n}$ the   It\^{o} measure on Brownian excursions and ${\mathfrak n}^+$  (resp. ${\mathfrak n}^-$)  the  restriction of  ${\mathfrak n}$ on  positive   excursions  (resp. negative   excursions).  For any  loop $\gamma$, let   $\ell_\gamma^0: = \lim_{\varepsilon\to0} \int_0^{T(\gamma)} 1_{\{0< \gamma(t)< \varepsilon\}} \d t$ be  its (total) local time at $0$.  

In the following lemma,  we identify a Poisson point process with its atoms.
\begin{lemma}\label{l:intensity} Let $\delta>0$.   
\begin{enumerate}[(i)]
\item \label{up}  The collection $\{(q,{\mathfrak e}_q^\uparrow),\, q\in {\mathcal Q}_{\delta \over 2}^\uparrow \}$   is a Poisson point process of intensity measure $\delta  {\rm d} a  \otimes   {\mathfrak n}^+(\d {\mathfrak e})$.
\item \label{down} The collection  $\{(q,{\mathfrak e}_q^\downarrow),\, q \in {\mathcal Q}^{\downarrow}_{\delta \over 2} \hbox{  such that }    q+ {\mathfrak e}_q^\downarrow \subset (0, \infty)  \}$   is a Poisson point process of intensity  measure $\delta  1_{\{a>0\}}{\rm d} a \otimes 1_{\{\min {\mathfrak e}> - a\}}   {\mathfrak n}^-(\d {\mathfrak e}) $. 
\item \label{up0}  The collection   $\{\min \gamma ,\,\gamma \in {\mathcal L}_{\delta\over 2} \hbox{ such that }0\in \gamma \}$   is a Poisson point process of intensity measure $ {\delta \over 2|a|}1_{\{a<0\}}  {\rm d} a  $.
\item  \label{loctime} Let $m>0$.  
 The collection  $ \{\ell_\gamma^0,\, \gamma \in {\mathcal L}_{\delta\over 2} \hbox{ such that } \min \gamma \in [-m,0],\, 0\in \gamma \}  $   is a Poisson point process of intensity measure $1_{\{\ell >0\}} {\delta \over 2 \ell} e^{-\ell / 2 m} {\rm d}\ell$. 
\item \label{loctime2} The collection  $\{ (\ell_\gamma^0, \, \gamma) ,\,\gamma \in {\mathcal L}_{\delta\over 2} \hbox{ such that }0\in \gamma \}$   is a Poisson point process of intensity measure $\frac\delta2\,  1_{\{\ell>0\}} \, {{\rm d} \ell \over \ell}  \p^*( (B_t,\, 0\le t\le \tau^B_\ell) \in {\rm d} \gamma) $, where   $\tau^B_\ell:= \inf\{s>0:  {\mathfrak L}_s> \ell\}$ denotes the  inverse of the Brownian local time,  and $\p^*( (B_t,\, 0\le t\le \tau^B_\ell) \in \bullet)$ is the projection of $\p( (B_t,\, 0\le t\le \tau^B_\ell) \in \bullet)$ on the space of unrooted loops. 
\end{enumerate}
\end{lemma}

\noindent {\it Proof}.  Item \eqref{up}  is Proposition 3.18, p. 44 of \cite{lupu}. Item \eqref{down} follows from the equality in distribution $\{(q,{\mathfrak e}_q^\downarrow),\, q \in {\mathcal Q}^{\downarrow}_{\frac{\delta}2}\}\law \{ - (q,{\mathfrak e}_q^\uparrow),\, q \in {\mathcal Q}^{\uparrow}_{\frac{\delta}2}\}$ and \eqref{up}.  Item \eqref{up0} comes from \eqref{up} and the fact that ${\mathfrak n}^+ (  r \in {\mathfrak e}) ={1\over 2 r}$ for any $r>0$ where, here and in the sequel, ${\mathfrak e}$ denotes a Brownian excursion. We prove now \eqref{loctime}.  The intensity measure  is given by
$$
\delta \int_{-m}^0 {\mathfrak n}^+ \big(  \ell^{|a|}_{\mathfrak e} \in {\rm d } \ell,\, |a| \in {\mathfrak e}  \big) \, \d  a
$$
where $\ell^r_{\mathfrak e}$ denotes the local time at $r$ of the excursion ${\mathfrak e}$.  Under ${\mathfrak n}^+$, conditionally on $|a|\in {\mathfrak e}$, the excursion after hitting $|a|$ is a Brownian motion killed at $0$. Therefore 
$$
{\mathfrak n}^+ \big(  \ell^{|a|}_{\mathfrak e} \in {\rm d } \ell \, \big|\, |a| \in {\mathfrak e}  \big) = \p_{|a|} ( {\mathfrak L}_{T_0^B}^{|a|} \in {\rm d}\ell ) = {1\over 2 |a| } e^{ - {\ell \over 2|a| }}  \, \d  \ell,
$$

\noindent where under $ \p_{|a|} $,  the Brownian motion $B$ starts at $|a|$ and ${\mathfrak L}_{T_0^B}^{|a|}$ denotes its  local time at  position  $|a|$   up to   $T_0^B:=\inf\{t>0: B_t=0\}$,  and the last equality follows from the standard Brownian excursion theory.   Hence the intensity measure is given by
$$
\delta \int_{-m}^0  {1\over  4a^2} e^{ - {\ell \over 2|a|} }  \, \d \ell\,  \d  a =  {\delta \over 2 \ell} e^{-\ell / 2 m}  \, \d \ell.   
$$
Finally, \eqref{loctime2} comes from Corollary 3.12, equation (3.3.5) p. 39 of \cite{lupu}.
$\Box$

\subsection{The Poisson--Dirichlet distribution}\label{s:poisson}

For a vector $\mathcal D=(D_1,D_2,\ldots)$ and a real $r$, we denote by $r {\mathcal D}$ the vector $(rD_1,rD_2,\ldots)$.  We recall  that for $a,b>0$, the density of the gamma$(a,b)$ distribution  is given by
$$
{1\over \Gamma(a)b^a} x^{a-1} e^{-{x\over b}}1_{\{x>0\}},
$$
and the density of the beta($a,b$) distribution is
$$
{\Gamma(a+b) \over \Gamma(a)\Gamma(b)}x^{a-1}(1-x)^{b-1}1_{\{x\in (0,1)\}}.
$$

We introduce the Poisson--Dirichlet distribution,  relying on Perman, Pitman and Yor \cite{ppy92}. Let $\beta>0$. Consider  a Poisson point process of intensity measure $ {\beta \over x} e^{-x}1_{\{x>0\}} \,{\rm d} x$ and denote by $\Delta_{(1)} \ge \Delta_{(2)} \ge \ldots $ its atoms. We can see them also as the jump sizes, ordered decreasingly, of a gamma subordinator of parameters $(\beta,1)$ up to time $1$. The sum $T:=\sum_{i\ge 1} \Delta_{(i)}$  has a gamma($\beta,1$) distribution. The random variable on the infinite simplex defined by 
 $$
 (P_{(1)}, P_{(2)},\ldots) := \left(\frac{\Delta_{(1)}}{T}, \frac{\Delta_{(2)}}{T}, \ldots\right)
 $$

\noindent has the Poisson--Dirichlet distribution with parameter $\beta$ (\cite{kingman}), and is independent of $T$ (\cite{mccloskey}, also Corollary 2.3 of \cite{ppy92}).  

Consider a decreasingly ordered positive vector  $(\xi_{(1)},\xi_{(2)},\ldots)$ of finite sum $\sum_{i\ge 1} \xi_{(i)}<\infty$. A size-biased random permutation, denoted by $(\xi_1,\xi_2,\ldots)$, is a permutation of $ (\xi_{(1)}, \xi_{(2)},\ldots)$  such that, conditionally on $\xi_1=\xi_{(i_1)},\ldots, \xi_j=\xi_{(i_j)}$, the term $\xi_{j+1}$ is chosen to be $\xi_{(k)}$ for $k\notin \{i_1,\ldots,i_j\}$ with probability ${ \xi_{(k)} \over \sum_{i\ge 1} \xi_{(i)} - (\xi_1+\ldots+\xi_j) }$.  The indices $(i_j,\, j\ge 1)$ can be constructed by taking i.i.d. exponential random variables of parameter $1$, denoted by $(\varepsilon_i,\, i\ge 1)$, and by ordering  $\mathbb{N}$   increasingly with respect to the total order $k_1 \leq k_2$ if and only if  $\xi_{(k_1)}/\varepsilon_{k_1} \ge \xi_{(k_2)}/\varepsilon_{k_2}$ (Lemma 4.4 of \cite{ppy92}). A result from McCloskey \cite{mccloskey} says that the $(P_i,\, i\ge 1)$ obtained from the $ (P_{(1)}, P_{(2)},\ldots)$ by size-biased ordering can also be obtained via the stick-breaking construction:
$$
P_i = (1-U_i)\prod_{j=1}^{i-1} U_j
$$
 where $U_i,\, i\ge 1$ are i.i.d. with law beta($\beta,1$). Let   \begin{equation}\label{def-Dbeta} {\mathcal D}_\beta:=(D_1,D_2,\ldots) \end{equation}
 
 \noindent be the point measure in $[0,1]$ defined by $D_i:=\prod_{j=1}^{i} U_j$, for $i\ge 1$.  
  
\begin{lemma}\label{l:D}
Let $m > 0 $, $\beta>0$ and $\Xi_\beta$ be a Poisson point process of intensity measure  ${ \beta \over a}1_{\{a>0\}} {\rm d} a  $. We denote by $a_1^{(m)}>a_2^{(m)}>\ldots$ the points of $\Xi_\beta$ belonging to $[0,m]$. Then $\big(a_i^{(m)},\, i\ge 1 \big)$ is distributed as $m {\mathcal D}_\beta$.
\end{lemma}

\noindent {\it Proof}.   For $0\le a \le 1$,    $$\p\big(\frac{1}{m} a_1^{(m)}   \le a\big) = \exp\big(- \int_{a m}^m \frac{\beta}{x} {\rm d} x\big)= a^\beta.$$

\noindent
Therefore it is a beta($\beta,1$) distribution.  Conditionally on  $\{a_1^{(m)}= a\}$, the law of $\Xi_\beta$ restricted to  $[0,a_1^{(m)})$ is the one of $\Xi_\beta $ restricted to $[0, a)$.  By iteration we get the Lemma.  $\Box$

\medskip
Denote by $L(t,r)$, $r\in \r$ and $t\ge 0$,   the local time of $X$ at  time $t$ and position $r $.   Let   \begin{equation}\label{def-tau}\tau_{r}(t):=\inf\{ s\ge 0\,:\, L(s,r)> t \} ,  \end{equation}

\noindent be the  inverse local time of $X$.   Denote by  \begin{equation}\label{def-T} T_{r}:=\inf\{ t\ge 0 \,:\, X_{t} = r \}
\end{equation} 

\noindent the hitting time of $r$.  We are interested in the process $(J(x))_{x\ge 0} $ defined as follows:  \begin{equation}\label{def-Jx} J(x) := \inf\{X_s,\, s\le \tau_0(x)\}, \qquad  x\ge 0. \end{equation}

Observe that under $\p^\delta$, $J$ is a Markov process. It has been studied in Section 4 of \cite{carmona-petit-yor94a}. The perturbed reflecting Brownian motions are related to the Poisson--Dirichlet via the following proposition.
\begin{proposition}\label{p:J} Let $\delta>0$. 
Under $\p^\delta$, the range  $\{J(x),\, x> 0\}$ is distributed as $-\Xi_\beta$ with $\beta={\delta\over 2}$.  Consequently, for any $m>0$, the range of $J$ in $[-m,0]$, ordered increasingly, is distributed as $-m {\mathcal D}_\beta$.
\end{proposition}

{\noindent\it Proof}.  Recall Proposition \ref{p:loopsoup}.  Note that  the range  $\{J(x),\, x> 0\}$  is equal to $\{ \min \gamma\, : \,0 \in \gamma, \, \gamma \in  {\mathcal L}_{\delta\over 2} \}$, hence the first statement is \eqref{up0} of Lemma \ref{l:intensity}. The second statement is Lemma \ref{l:D}. $\Box$

\medskip

Under $\p^\delta$,  for $m>0$, let \begin{equation}\label{def-TJ} T_{-m}^J:= \inf\{x>0: J(x) \le -m\} \end{equation}

\noindent be the first passage time of $-m$ by $J$.  The main result in this subsection is the following description of the jump times of $J$ before its first passage time of $-m$:

\begin{proposition}\label{p:Jjump}
For $m>0$, let   $x_1^{(m)}> x_2^{(m)}>\ldots$ denote the jumping times of $J$  before time $T_{-m}^J$. Under $\p^\delta$: 
\begin{enumerate}[(i)]
\item (\cite{yor92})  \label{gamma}$T_{-m}^J$ follows a gamma(${\delta \over 2},2m$) distribution. Consequently, for any $x>0$, ${-1\over J(x)}$ follows a gamma(${\delta \over 2}, {2\over x}$) distribution.
\item $T_{-m}^J$ is independent of $ {1\over T_{-m}^J} (x_1^{(m)},x_2^{(m)},\ldots )$.
\item $ {1\over T_{-m}^J} (x_1^{(m)},x_2^{(m)},\ldots )$ is distributed as ${\mathcal D}_\beta$ with $\beta={\delta \over 2}$.
\end{enumerate}
\end{proposition}

Statement \eqref{gamma} of Proposition \ref{p:Jjump}  is not new. It is contained in Proposition 9.1, Chapter 9.2, p. 123, of Yor \cite{yor92}.  For the sake of completeness we give here another proof of \eqref{gamma} based on Lemma \ref{l:intensity}.

\medskip

\noindent {\it Proof}.   (i)    
    Lemma \ref{l:intensity}, \eqref{loctime}   says that  $ \{\ell_\gamma^0,\, \gamma \in {\mathcal L}_{\delta\over 2} \hbox{ such that } \min \gamma \in (-m,0],\, 0\in \gamma \}  $   forms  a  Poisson point process  of intensity measure   $1_{\{\ell >0\}} {\delta \over 2 \ell} e^{-\ell / 2 m} {\rm d}\ell$,     whose atoms are exactly   the (non-ordered) sequence  $\{x_{i-1}^{(m)}-x_i^{(m)}, i\ge 1\}$ where $x_0^{(m)}:= T_{-m}^J$.  Note that  $T^J_{-m}= L(T_{-m}, 0)= \sum_{\min \gamma \in (-m,0],\, 0\in \gamma} \ell_\gamma^0$.

 Let for $i\ge 1$, $d_i^{(m)}:=(x_{i-1}^{(m)}-x_i^{(m)})/T_{-m}^J$ and denote by  $ \{d_{(1)}^{(m)}>d_{(2)}^{(m)}>\ldots \}$  the sequence ordered decreasingly.   
 Then the properties of the Poisson--Dirichlet distribution recalled at the beginning of the section imply that $T_{-m}^J$ is independent of the point measure $ \{d_{(1)}^{(m)},d_{(2)}^{(m)},\ldots \}$  and that $T_{-m}^J/2m$ follows a gamma(${\delta \over 2},1$) distribution. Also, the second statement of \eqref{gamma} comes from the observation that $\{ J(x) > -m \} = \{T_{-m}^J > x \}$. This proves  \eqref{gamma}.

 (ii) and (iii): 
  It remains to show that the vector $(d_1^{(m)},d_2^{(m)},\ldots)$  is a size-biased ordering of  $ \{d_{(1)}^{(m)},d_{(2)}^{(m)},\ldots \}$, and that this size-biased ordering is still independent of $T_{-m}^J$.

 To this end, denote by $\{(-m_i,\, i\in  \mathcal{I})\}$ the range of $J$.   By Proposition  \ref{p:J},  the point measure $\{ \ln(m_i),\, i\in  \mathcal{I}\}$ is a Poisson point process on $\r$ of intensity measure ${\delta \over 2} \d t$.  
 
 When $J$ jumps at some $-m_i$, the time to jump at $-m_{i+1}$ is exponentially distributed with parameter ${\mathfrak n}^-( \min {\mathfrak e} \le - m_i) ={1\over 2 m_i}$ (it is the local time at $0$ of a Brownian motion when it hits level $-m_i$, by Markov property of the process $(X,I)$ under $\p^\delta$). 
  
 Denote by $\varepsilon_i$ the exponential of parameter $1$  obtained as the waiting time between jumps to $-m_i$ and to $-m_{i+1}$, divided by $2 m_i$. Conditionally on $\{(m_i,\, i\in  \mathcal{I})\}$, the random variables ${(\varepsilon_i,\, i\in  \mathcal{I})}$ are i.i.d. and  exponentially distributed with parameter $1$.  Then $\{ (\ln(m_i),\varepsilon_i),\, i\in \mathcal{I}\}$ is a Poisson point process on $\r\times \r_+$ of intensity measure ${\delta\over 2} {\rm d} t \otimes e^{-x}{\rm d} x$.   It is straightforward   to  check that $\{ (\ln(2m_i\varepsilon_i),\varepsilon_i),\, i\in \mathcal{I}\}$ is still a  Poisson point process with the same intensity measure.

 Suppose that we enumerated the range of $J$ with $\mathcal{I}=\mathbb{Z}$ so that $(-m_i,\,i\ge 1)$ are the atoms of the range in $(-m,0)$ ranked increasingly. Then, $2m_i\varepsilon_i =  T_{-m}^{J} d_i^{(m)}=:\xi_i$ for any $i\ge 1$. We deduce that, conditionally on $\{T_{-m}^{J} d_i^{(m)}, i\ge 1\}$, the vector $(\varepsilon_1,\varepsilon_2,\ldots)$ consists of i.i.d.  random variables exponentially distributed with parameter $1$, and $i\le j$ if and only if $\xi_i/\varepsilon_i \ge \xi_j/\varepsilon_j$. From the description of size-biased ordering at the beginning of this section, we conclude that the $(\xi_i,\, i\ge 1)$ are indeed size-biased ordered, hence also $ (d_1^{(m)},d_2^{(m)},\ldots )$. $\Box$

 \medskip

  Since $T_{-m}^J= L(T_{-m},0)$, the statement  \eqref{gamma} of Proposition \ref{p:Jjump}   says     \begin{equation} \label{yor92} L(T_{-m}, 0)    \law {\mbox gamma}(\frac\delta2, 2 m ). \end{equation}

\begin{corollary}\label{c:jumptimes} Let $\delta>0$. 
Under $\p^\delta$, the collection of jumping times of $J$ is distributed as $\Xi_\beta$ with $\beta={\delta\over 2}$.
\end{corollary}

{\noindent\it Proof}.  From  Proposition \ref{p:Jjump} and Lemma \ref{l:D}, we can couple the jumping times of $J$ which are strictly smaller than the passage time of $-m$ with $\Xi_\beta$ restricted to $[0, Z_m]$ where $Z_m$ is  gamma(${\delta \over 2}, 2m$) distributed, independent of $\Xi_\beta$.  Letting  $m\to +\infty$ gives the Corollary. $\Box$

\subsection{Some known results}\label{s:known}

At first we recall two Ray--Knight theorems:

\begin{proposition}[Le Gall and Yor \cite{legall-yor}]\label{p:legall-yor}
Let $\delta>0$. Under $\p^{(-\delta)}$, the process $\big( L(\infty, t),\, t\ge 0 \big)$ is  the square of a Bessel process of dimension $\delta$  starting from $0$ reflected at $0$.
\end{proposition}

 \begin{proposition} [Carmona, Petit and Yor \cite{carmona-petit-yor94b},  Werner \cite{werner95}]  \label{p:carmona-petit-yor}
Let $\delta>0$ and $a\ge 0$. Under $\p^\delta$, the process $\big( L(\tau_{0}(a), -t),\, t\ge 0 \big)$ is   the square of a Bessel process of dimension $(2-\delta)$ starting from $a$ absorbed at $0$. 
 \end{proposition}

 We have the following independence result.

\begin{proposition}[Yor \cite{yor92}, Proposition 9.1] \label{p:min}   Let $\delta>0$. Under $\p^\delta$, 
for any fixed $x>0$, $L(T_{J(x)}, 0)/x$ is  independent of $J(x)$ and follows a beta(${\delta\over 2},1$) distribution. 
\end{proposition}

 We introduce some notations which will be used in Section \ref{s:decompositionminimum}. 

\begin{notation}\label{n:<h}
Let $h \in \r$.  We define the process $X^{-,  h}$ obtained by gluing the excursions of $X$ below $h$ as follows.    Let for $t\ge0$, 
$$
A_t^{ -,  h} := \int_0^t 1_{\{X_s \le  h \}} {\rm d} s,\, \qquad  \alpha_t^{ -,  h} := \inf\{u>0,\, A_u^{-, h} > t\},$$
with the usual convention   $\inf \emptyset:=\infty$.  Define $$ X^{-,  h}_t := X_{\alpha_t^{-,  h}}, \qquad   t < A_\infty^{-,  h}:=\int_0^\infty 1_{\{X_s \le  h\} } {\rm d} s.$$ 

Similarly, we define   $A_t^{+,  h}$, $\alpha_t^{+,  h}$ and $X^{+,  h}$ by replacing $X_s \le h$ by $X_s > h$. When the process is denoted by $X$ with some superscript, the analogous quantities will  hold the same superscript. For example for $r\in \r$,  $\ell >0$,  $\tau_r^{+,  h}(\ell)=\inf\{t>0: L^{+, h}(t, r) > \ell\}$, where $L^{+, h}(t, r)$ denotes the  local time of $X^{+,  h}$  at position $r$ and time $t$.  
\end{notation}

\begin{proposition}[Perman and Werner \cite{perman-werner}] \label{p:perman-werner} Let $\delta>0$. 
Under $\p^\delta$, the two processes $X^{+, 0}$ and $X^{-, 0}$ are independent. Moreover, $X^{+,  0}$ is a reflecting Brownian motion, and the process $(X_t^{-,  0},\inf_{s\le t} X_s^{-,  0} )_{t\ge 0}$ is strongly Markovian.
\end{proposition}

Let $m>0$. We look at these processes up to the first time the process $X$ hits level $-m$. In this case, there is a   dependence between $X^{-, 0}$ and $X^{+,  0}$ due to their duration. This dependence is taken care of by conditioning on the (common) local time at $0$ of $X^{-, 0}$ and $X^{+,  0}$. It is the content of the following corollary.
\begin{corollary}\label{c:ind} Let $\delta>0$.  
Fix $m>0$. Under $\p^\delta$, conditionally on $(X^{-, 0}_t,\,t\le A_{T_{-m}}^{-,  0}  )$, the process  $(X^{+, 0}_t, t\le  A_{T_{-m}}^{+, 0})$  is a reflecting Brownian motion stopped at time $\tau_0^{+, 0}(\ell)$ where  $\ell=L^{-, 0}(A^{-,  0}_{T_{-m}},0)=L(T_{-m}, 0)$.
\end{corollary}

{\noindent\it Proof}.  By    Proposition \ref{p:perman-werner}, conditionally on $X^{-, 0}$, the process $(X^{+, 0}_{\tau_0^{+, 0}(t)},\, t\ge 0)$ is a reflecting Brownian motion indexed by its inverse local time. Observe that $A_{T_{-m}}^{+, 0} = \tau_0^{+, 0}(\ell)$ with $\ell = L^{-, 0}(0,T_{-m}^{-, 0})$.  It proves the Corollary. $\Box$

 \medskip

As mentioned in Section 3 of Werner \cite{werner95}, we have the following duality between $\p^{(-\delta)}$ and $\p^\delta$.
\begin{proposition}[Werner \cite{werner95}] \label{p:duality} Let $\delta>0$. 
For any $m>0$, the process $(X_{T_{-m}-t}+m,\, t\le T_{-m})$ under $\p^\delta$ has the distribution of $(X_t,\, t \le {\mathscr D}_m  )$ under $\p^{(-\delta)}$, where ${\mathscr D}_m:= \sup\{t>0: X_t=m\}$ denotes the last passage time at $m$.
\end{proposition}

\section{ Decomposition at a hitting time}

\label{s:decomposition}

  The following lemma is  Lemma 2.3 in Perman \cite{perman}, together with the duality stated in Proposition  \ref{p:duality}.  Recall  that  under $\p^{(-2)}$, $X$ is a Bessel process of dimension 3.
  We refer to  \eqref{def-T} for the definition of the first hitting time $T_r$ and to \eqref{def-Jx} for the process $J(x)$. 
  
 \begin{lemma}[Perman \cite{perman}] \label{l:perman} Let $\delta>0$. 
 Let $m,x>0$ and $y\in (0,x)$.  Define the  processes 
$$
Z^1 := (X_{T_{J(x)} - t} - J(x))_{t\in [0,T_{J(x)} ]} \qquad Z^2:= (X_{T_{J(x)} + t} - J(x))_{t\in [0,\tau_0(x) - T_{J(x)} ]}.
$$
Under $\p^\delta (\cdot | J(x)=-m, L(T_{J(x)}, 0) =y)$: 
\begin{enumerate}[(i)]
\item $Z^1$ and $Z^2$ are   independent, 
\item $Z^1$ is distributed as $(X_t)_{t \in [0, {\mathscr D}_m]} $ under $\p^{(-\delta)} ( \cdot | L(\infty, m) = y )$,
\item $Z^2$ is distributed as $(X_t)_{t \in [0, {\mathscr D}_m]}$ under $\p^{(-2)}( \cdot | L(\infty, m) = x-y )$,
\end{enumerate}
\noindent with  ${\mathscr D}_m:= \sup\{t>0: X_t=m\}$.
 \end{lemma}

{\noindent\it Proof.}  For the sake of completeness, we give here a proof which is different from  Perman \cite{perman}'s.  

 By Proposition \ref{p:loopsoup}, we can identify under $\p^\delta$ the point measure $\{ (q, {\mathfrak e}_{X,q}^\uparrow),\, q \in {\mathcal Q}^\uparrow_X\}$ with $\{ (q, {\mathfrak e}_q^\uparrow),\, q \in {\mathcal Q}^\uparrow_{\delta \over 2} \cap (-\infty,0)\}$.  Using the notations in Lemma \ref{l:intensity},  we  have $$ 
L(T_{J(x)}, 0) 
=
\sum_{ q  \in {\mathcal Q}^\uparrow_{\delta \over 2} \cap  ( J(x) , 0) } \ell^0_{\gamma} \,  <   x  \, \le  \sum_{q  \in {\mathcal Q}^\uparrow_{\delta \over 2} \cap   [J(x) , 0) } \ell^0_{\gamma} , 
$$

 \noindent where  in the above sum $\gamma$  is the (unique) loop in $\mathcal {L}_{{\delta\over 2}}$ such that $\min \gamma =q$.  Let $\ell^r({\mathfrak e})$ be the local time of the excursion ${\mathfrak e}$ at level $r$.  We claim that conditioning  on $\{J(x) = -m,  L(T_{J(x)},0) =y\}$, $   {\mathfrak e}_{J(x)}^\uparrow  $  and $\{ (q, {\mathfrak e}_q^\uparrow),\, q  \in {\mathcal Q}^\uparrow_{\delta \over 2} \cap (J(x), 0)\}$ are independent and distributed  as  a Brownian excursion  ${\mathfrak e}$  under ${\mathfrak n}^+ ( \cdot \, |\, \ell^{m}({\mathfrak e}) > x- y)$, and $\{ (q, {\mathfrak e}_q^\uparrow),\, q  \in  {\mathcal Q}^\uparrow_{\delta \over 2} \cap (-m, 0)\}$ conditioned on $\{\xi_m= y\}$ respectively, where $\xi_m:=\sum_{q \in {\mathcal Q}^\uparrow_{\delta \over 2} \cap (-m, 0)} \ell^0_{\gamma}$. 

In fact,    let $F: \r_- \times {\mathcal K}  \to \r_+ $,   $G: {\mathcal K}  \to \r_+$ and $ f: \r^2 \to \r_+$ be three measurable functions.     Note that 
$ \xi_m = \sum_{q\in {\mathcal Q}^\uparrow_{\delta \over 2} \cap (-m, 0)} \ell^{|q|} ({\mathfrak e}_q^\uparrow) .$  By   Proposition \ref{p:loopsoup},  we deduce from the master formula that  \begin{eqnarray*}
&&
\e^\delta \left[  \ee^{-\sum_{q \in {\mathcal Q}^\uparrow_{\delta \over 2} \cap (J(x), 0)}F(q, {\mathfrak e}_q^\uparrow)} \, G({\mathfrak e}_{J(x)}^\uparrow) \, f(J(x), L(T_{J(x)}, 0) ) \right]
\\
&=&
\e^\delta \left[\sum_{m>0} \ee^{-\sum_{q \in {\mathcal Q}^\uparrow_{\delta \over 2} \cap (-m, 0)}F(q, {\mathfrak e}_q^\uparrow)}  \, G({\mathfrak e}_{-m}^\uparrow) \, f(-m, \xi_m) \, 1_{\{\xi_m < x,  \ell^m({\mathfrak e}_{-m}^\uparrow)  > x- \xi_m\}}\right]
\\
&=&
\delta\, \int_0^\infty \d m\,  \e\left[\ee^{-\sum_{q \in {\mathcal Q}^\uparrow_{\delta \over 2} \cap (-m, 0)}F(q, {\mathfrak e}_q^\uparrow)}  \,   f(-m, \xi_m)\, 1_{\{\xi_m < x\}} \, \int {\mathfrak n}^+ (\d {\mathfrak e}) G({\mathfrak e}) \, 1_{\{\ell^m({\mathfrak e})  > x- \xi_m\}}\right],
\end{eqnarray*}

\noindent by using Lemma \ref{l:intensity} (i).  The claim follows.

  Now we observe  that  $Z^2$ is measurable with respect to $ {\mathfrak e}_{J(x)}^\uparrow $ whereas $Z^1$ is to $\{ (q, {\mathfrak e}_q^\uparrow),\, q  \in {\mathcal Q}^\uparrow_{\delta \over 2} \cap (J(x), 0)\}$. It yields  (i). Moreover,  conditioning  on $\{J(x) = -m,  L(T_{J(x)}, 0) =y\}$, $Z^2$ is distributed as $({\mathfrak e}_t)_{t\in [0, \sigma^m_{x-y}]}$, under ${\mathfrak n}^+ ( \cdot \, |\, \ell^{m}({\mathfrak e}) > x- y)$, where $\sigma^m_{x-y}:=\inf\{t>0: \ell^m_t({\mathfrak e})=x-y\}$ with $\ell^m_t({\mathfrak e})$ being the local time at level $m$ at time $t$.  The latter process  has the same law as $(X_t)_{t \in [0, {\mathscr D}_m]}$ under $\p^{(-2)}( \cdot | L(\infty, m) = x-y )$.\footnote{Under ${\mathfrak n}^+ ( \cdot \, |\, \ell^{m}({\mathfrak e}) > x- y)$, an excursion  up to the inverse local time $x-y$ at position $m$  is a three-dimensional Bessel process, up to the hitting time of $m$, followed by a Brownian motion starting at $m$ stopped at local time at level $m$ given by $x-y$, this Brownian motion being conditioned on not touching $0$ during that time. By excursion theory, the  time-reversed process is distributed as a Brownian motion starting at level $m$ stopped at the hitting time of $0$ conditioned on the local time at $m$ being equal to $x-y$. We conclude by William's time reversal theorem   (Corollary VII.4.6 of \cite{revuz-yor}). } We get (iii).

To prove (ii), we denote by $\widehat e$ the time-reversal of a loop $e$. 
By Proposition \ref{p:duality},   $\{ (m+q,\widehat{\mathfrak e}_{q,X}^\uparrow),\, q  \in {\mathcal Q}^\uparrow_X\cap (-m, 0)\}$ under $\p^\delta$ is distributed as $\{ (q, {\mathfrak e}_{q,X}^\uparrow),\, q  \in {\mathcal Q}^\uparrow_X\cap (0,m)\}$ under $\p^{(-\delta)}$. Note that under $\p^\delta(\cdot | J(x) = -m,  L(T_{J(x)}, 0) =y)$,  $Z^1$ can be constructed from $\{ (m+q, {\widehat {\mathfrak e}_{q,X}^\uparrow)},\, q  \in {\mathcal Q}^\uparrow_X\cap (-m, 0)\}$. Then $Z^1$ is distributed as $(X_t)_{0\le t \le {\mathscr D}_m}$ under  $\p^{(-\delta)}( \cdot |  \sum_{q\in {\mathcal Q}^\uparrow_X\cap(0, m)}  \ell^{m-q} ({\mathfrak e}_{q,X}^\uparrow)=y  )$.  Finally remark that $ \sum_{q\in {\mathcal Q}^\uparrow_X\cap(0, m)} \ell^{m-q} ({\mathfrak e}_{q,X}^\uparrow) =L({\mathscr D}_m, m)= L(\infty, m)$. We get (ii).  This completes the proof of Lemma \ref{l:perman}. $\Box$

 \medskip

Fix $m>0$. 
The following Theorems \ref{p:hitting} and \ref{p:hitting2} describe the path decomposition of $(X_t)$ at $T_{-m}$.      Let ${\tt g}_m := \sup\{t\in[0,T_{-m}]  \, : \,  X_t  = 0   \}$.   Recall that  $I_{{\tt g}_m}= \inf_{0\le s \le {\tt g}_m} X_s$.  Define $${\tt  d}_m:= \inf\{t> {\tt g}_m : X_t= I_{{\tt g}_m}\}.$$

\begin{theorem}\label{p:hitting} Let $\delta>0$. 
Fix $m>0$. Under $\p^\delta$, the random variable $\frac1{m} |I_{{\tt g}_m}|$ is beta$({\delta \over 2},1)$ distributed.   Moreover, for $0<a<m$, conditionally on $ \{I_{{\tt g}_m}=-a\}$,  the four processes  
\begin{align*}
& (X_t,\, t\in [0,T_{-a}]) ,\\   &(X_{{\tt g}_m-t},\, t\in [0,{\tt g}_m-T_{-a}] , \\  &(- X_{{\tt g}_m + t},\, t\in [0 , {\tt  d}_m  - {\tt g}_m ]) ,\\ &(X_{{\tt d}_m + t} + a,\, t\in [0 , T_{-m}  - {\tt d}_m  ]) ,
\end{align*}
are independent, with law respectively the one of:
\begin{enumerate}[(i)]
\item  $X$ under $\p^\delta$ up to the hitting time of $-a$;
\item a Brownian motion up to the hitting time of $-a$;
\item  a Bessel process of dimension 3 from $0$ stopped when hitting $a$;
\item   $X$ under $\p^\delta$ conditionally on $\{T_{-(m-a)} < T_a\}$. 
\end{enumerate}
\end{theorem}

{\noindent\it Proof}.  To get the distribution of  $\frac1{m} I_{{\tt g}_m}$  we proceed as follows:    under $\p^\delta$, $X$ is measurable with respect to  its excursions above the infimum, that we denoted by $({\mathfrak e}_{q,X}^\uparrow,\, q\in {\mathcal Q}_X^\uparrow)$, that we identify with $({\mathfrak e}_q^\uparrow,\, q\in {\mathcal Q}_{\delta\over 2}^\uparrow)$ by Proposition \ref{p:loopsoup}.   The variable $ I_{{\tt g}_m}$ is the global minimum of the loops $\gamma$ such that $\min \gamma> -m $ and $0\in \gamma$.  By Lemma \ref{l:intensity} (iii), we get the law of $ I_{{\tt g}_m}$ (it is  also a consequence of Lemma \ref{l:D}  together with Proposition \ref{p:J}).

Let $\tilde \gamma$ be the loop such that $\min \tilde \gamma =  I_{{\tt g}_m}$ and call $-a= I_{{\tt g}_m}$ its minimum. Conditioning on $\tilde \gamma$ and loops hitting $(-\infty,-a)$,  the loops $\gamma$ such that $\min \gamma > -a$ are distributed as the usual Brownian loop soup $\mathcal L_{\delta \over 2}$ in $(-a,\infty)$. It gives (i) by Proposition \ref{p:loopsoup}.  Conditioning on $\min \tilde \gamma= I_{{\tt g}_m}=-a$ and on loops hitting $(-\infty,-a)$, the loop $\tilde \gamma - \min \tilde \gamma$ has the measure ${\mathfrak n}^+( {\rm d} {\mathfrak e}| \max {\mathfrak e} > a )$. Therefore (ii) and (iii) come from the usual decomposition of the It\^o measure.  Finally conditioning on $\min \tilde \gamma= I_{{\tt g}_m}=-a$, the collection of loops  $\gamma$ with $\min \gamma \in (-m,-a)$ is  distributed as the Brownian loop soup $\mathcal L_{\delta \over 2}$ restricted to loops $\gamma$ such that $\min \gamma \in (-m,-a)$ conditioned on the event that none of these loops hit $0$. We deduce (iv). $\Box$

\medskip

The following theorem gives the path decomposition  when conditioning on $(L(T_{-m}, 0),  I_{{\tt g}_m})$. Recall \eqref{yor92} for the law of $L(T_{-m}, 0)$.

\begin{theorem}\label{p:hitting2}  We keep the notations of  Theorem \ref{p:hitting}.  Under $\p^\delta$,  

(i)   the   density of $(L(T_{-m}, 0), |I_{{\tt g}_m}|)$ is given by $ \frac{a^{-2}}{2 \Gamma(\frac\delta2)}  (\frac{x}{2m})^{\frac\delta2} e^{- \frac{x}{2a}} $ for $x>0$ and $0< a <m$.

(ii)  conditionally on $\{L(T_{-m}, 0)=x, I_{{\tt g}_m}= -a \}$, the three processes
\begin{align*}
&   (X_t , \,  t \in [ 0 , {\tt g}_m ]), \\   &   (- X_{{\tt g}_m+t},\, t\in [0, {\tt d}_m - {\tt g}_m]),  \\  &    (X_{{\tt d}_m + t} +a ,\, t \in [0, T_{-m}-{\tt d}_m]),
\end{align*}

\noindent are independent and distributed respectively as 
\begin{eqnarray*}
&& (X_t , t\in [0,  \tau_0(x)]) \mbox{ under } \p^\delta(\cdot \, |\, J(x)=-a),
\\
&&  \mbox{a Bessel process of dimension 3 starting from $0$ stopped when hitting $a$}, \\
&& \mbox{$X$ under $\p^\delta$ conditionally on $\{T_{-(m-a)} < T_{a}\}$}; 
\end{eqnarray*}

\end{theorem}

{\noindent\it Proof.}     By Theorem \ref{p:hitting}, conditionally on $\{ I_{{\tt g}_m}= -a \}$,  the three processes
\begin{align*}
&   (X_t , \,  t \in [ 0 , {\tt g}_m ]), \\   &   (- X_{{\tt g}_m+t},\, t\in [0, {\tt d}_m - {\tt g}_m]),  \\  &    (X_{{\tt d}_m + t} +a ,\, t \in [0, T_{-m}-{\tt d}_m]),
\end{align*}

\noindent are independent. Since $L(T_{-m},0)$ is measurable with respect to $\sigma (X_t , \,  t \in [ 0 , {\tt g}_m ])$, we obtain the independence of the three processes in (ii) and the claimed  laws of the last two processes in (ii).

To complete the proof,      it is enough  to show  that  for any    bounded continuous functional $\Phi$ on ${\mathcal K}$ and  any bounded continuous  function $f: \r^2 \to \r $,  \begin{eqnarray} && \e^\delta [\Phi(X_t , \,  t \in [ 0 , {\tt g}_m ]) f(L(T_{-m}, 0),  I_{{\tt g}_m})] 
\label{condLI}  
\\
&=&
\int_0^\infty \int_0^m     \e^\delta [ \Phi(X_t , t\in [0,  \tau_0(x)]) \, |\, J(x)= -a]  f(x,  -a)    \frac{a^{-2}}{2 \Gamma(\frac\delta2)}  (\frac{x}{2m})^{\frac\delta2} e^{- \frac{x}{2a}} \d a \d x . \nonumber
\end{eqnarray}


   By  Theorem \ref{p:hitting},  \begin{eqnarray} && \e^\delta [\Phi(X_t , \,  t \in [ 0 , {\tt g}_m ]) f(L(T_{-m}, 0),  I_{{\tt g}_m})]
\nonumber 
\\
&=&
\int_0^m  \frac\delta2 m^{-\frac\delta2} a^{\frac\delta2-1}   \e^\delta [ \Phi(X^{1,a} \oplus X^{2,a}) f(L^0(X^{1,a}) +L^0(X^{2,a}),  -a)]  \d a  , \label{hitting1}
\end{eqnarray}


 \noindent where $X^{1,a}_s:= X_s, s \le T_{-a}$,  $X^{2, a}$ is the time-reversal of an independent Brownian motion up to its hitting time of $-a$  (so $X^{2, a}$ starts from  $-a$ and ends at $0$),  $X^{1,a} \oplus X^{2,a}$ denotes the process obtained by gluing $X^{2, a}$ and $X^{1,a}$ at time $T_{-a}$, and $L^0(X^{1,a})$ (resp. $ L^0(X^{2,a})$) is the local time  at position  $0$ of $X^{1,a}$ (resp.  $X^{2,a}$).

  The standard excursion theory  says that   $  \p^\delta (L^0(X^{2,a}) \in \d z)= \frac1{2a} e^{- \frac{z}{2a}}\d z,   z >0.$   By \eqref{yor92}, $L^0(X^{1,a}) \law {\rm gamma}(\frac\delta2, 2 a)$. 
Then for any    bounded Borel function $h$, we have    \begin{eqnarray*}
&& \e^\delta [ \Phi(X^{1,a}\oplus X^{2,a})  \, h( L^0(X^{1,a}) +L^0(X^{2,a}) ]
 \\
&=&
\int_0^\infty\int_0^\infty  h(y+z) \e^\delta [\Phi(X^{1,a}\oplus X^{2,a}) | L^0(X^{1,a})  = y, L^0(X^{2,a})= z]  \frac{(2 a)^{-1-\frac\delta2} }{\Gamma(\frac\delta2)} y^{\frac\delta2-1} e^{- \frac{y+z}{2a}} \d y \d z
\\
&=&
\int_0^\infty h(x) \int_0^x \e^\delta [\Phi(X^{1,a}\oplus X^{2,a}) | L^0(X^{1,a})  = y, L^0(X^{2,a})= x-y]  \frac{(2 a)^{-1-\frac\delta2} }{\Gamma(\frac\delta2)} y^{\frac\delta2-1} e^{- \frac{x}{2a}} \d y \d x
\\
&=&
\int_0^\infty h(x) \int_0^x \e^\delta [ \Phi(X_t, t \le \tau_0(x)) | J(x)= - a, L(T_{J(x)}, 0)=y]   \frac{(2 a)^{-1-\frac\delta2} }{\Gamma(\frac\delta2)}  y^{\frac\delta2-1} e^{- \frac{x}{2a}} \d y \d x ,
\end{eqnarray*}

\noindent where the last equality is due to   Lemma \ref{l:perman}.  Since $\p^\delta( L(T_{J(x)}, 0) \in \d y ) =  \frac\delta2 x^{-\frac\delta2}  y^{\frac\delta2-1}  1_{\{0< y < x\}} \d y$ (see Proposition \ref{p:min}),  we get that \begin{eqnarray*}
&& \e^\delta [ \Phi(X^{1,a}\oplus X^{2,a})  \, h( L^0(X^{1,a}) +L^0(X^{2,a}) ]
\\
&=&
\int_0^\infty h(x)   \e^\delta [ \Phi(X_t, t \le \tau_0(x)) | J(x)= - a]   \frac{(2 a)^{-1-\frac\delta2} }{\Gamma(1+\frac\delta2)} \, x^{\frac\delta2}\,   e^{- \frac{x}{2a}}   \d x ,
\end{eqnarray*}

\noindent which in view of \eqref{hitting1}   yields    \eqref{condLI} and  completes the proof of the Proposition. 
$\Box$

 \begin{remark} We may also directly prove (i) as follows: In view of \eqref{yor92}, it is enough to show \begin{equation} \label{l:Jcond} \p^\delta\big( | I_{{\tt g}_m}| \in \d a \,|\,  L(T_{-m}, 0)=x\big)   = \frac{x}{2 a^2} 
e^{-\frac{x}{2a}  + \frac{x}{2 m}} 1_{\{0< a < m\}} \d a.
\end{equation}

To this end, we  shall   prove that conditionally on $\{L(T_{-m}, 0)=x\}$, $ I_{{\tt g}_m}$ is distributed as $\inf_{0\le t \le  \tau_0^B(x)} B(t) $ conditioned on $\{\inf_{0\le t \le  \tau_0^B(x)} B(t)  > -m\}$, where $\tau_0^B(x):=\inf\{t>0: {\mathfrak L}_t> x\}$ denotes   the first time when the local time at $0$ of $B$ attains $x$. 

Consider the Brownian loop soup ${\mathcal L}_{\delta\over 2}$. In this setting (recalling Proposition \ref{p:loopsoup}), $$ I_{{\tt g}_m}
=
\inf_{\gamma \in {\mathcal L}_{\delta\over 2}} \Big\{ q :  q= \min\gamma > - m, 0 \in \gamma\Big\},$$

\noindent and $$L(T_{-m}, 0)
=
\sum_{\gamma \in {\mathcal L}_{\delta\over 2}}   \ell_\gamma^0 \,  1_{\{ \min \gamma \in (-m, 0), 0 \in \gamma \}}.$$

  From \eqref{loctime2} of Lemma \ref{l:intensity}, conditionally on $\{\ell_\gamma^0:  \gamma \in {\mathcal L}_{\delta\over 2}, 0 \in \gamma\} $, the loops $\gamma$ such that $0\in \gamma$ are (the projection on the space of unrooted loops of) independent Brownian motions stopped at $\tau_0^\ell$ with $\ell= \ell_\gamma^0$. Then,   the loops $\gamma$ such that $0\in \gamma$ and $\min\gamma>-m$ are merely (the projection of) independent Brownian motions stopped at local time given by $\ell_\gamma^0$, conditioned on not hitting $-m$. 
   The conditional density \eqref{l:Jcond} of   $ I_{{\tt g}_m}$   follows from standard Brownian excursion theory.  $\Box$
 \end{remark}
 
\begin{corollary} \label{c:hitting2} Let us keep the notations of  Theorem \ref{p:hitting}. 
Let $x> 0$ and $m>a >0$. Under $\p^\delta$,    the conditional law of  the process $  (X_t , \,  t \in [ 0 , {\tt g}_m ]) $  given $\{L(T_{-m}, 0)=x\}$ is equal to the (unconditional) law of  $(X_t, t\in [0, \tau_0(x)])$     biased by $c_{m, x, \delta} |J(x)|^{ {\delta \over 2} -1 } {1}_{\{J(x)> -m \}}$, with $$c_{m, x, \delta}:= \Gamma(\frac\delta2) \,  (\frac{x}{2})^{1- \frac\delta2} e^{  \frac{x} {2 m}} . $$
\end{corollary}

 {\noindent\it Proof.}  Let  $\Phi$ be a bounded continuous functional   on ${\mathcal K}$.     Recall from \eqref{yor92} that the density function of  $L(T_{-m}, 0)$ is $\frac1{\Gamma(\frac\delta2)} (2 m)^{-\frac\delta2} x^{\frac\delta2-1} e^{- \frac{x}{2m}}, x >0$. 
Considering some $f$ in \eqref{condLI} which only depends on the first coordinate, we see that for all $x>0$, \begin{eqnarray} && \e^\delta [\Phi(X_t , \,  t \in [ 0 , {\tt g}_m ]) \, |\, L(T_{-m}, 0) =x] 
\nonumber 
\\
&=&
\int _0^m      \e^\delta [ \Phi(X_t , t\in [0,  \tau_0(x)]) \, |\, J(x)= -a]     \frac{x}{2} a^{-2} e^{-  \frac{x}{2a} +\frac{x}{2m}}  \d a 
\nonumber \\
&=&
c_{m, x, \delta}\, \e^\delta [ \Phi(X_t , t\in [0,  \tau_0(x)]) \, |J(x)|^{\frac\delta2-1} \, 1_{\{J(x)>-m\}}]  ,  \label{eq:Phigm}
\end{eqnarray}

\noindent by using the fact that   the density  of $|J(x)|$ is  $a \to \frac1{\Gamma(\frac\delta2)} (\frac{x}{2})^{\frac\delta2} a^{- \frac\delta2-1} e^{- \frac{x} {2 a}}$.  This proves   Corollary \ref{c:hitting2}. $\Box$

\begin{remark}\label{r:Phigm}  Note that   the conditional expectation term on the left-hand-side of \eqref{eq:Phigm} is a continuous function of $(m, x)$, this fact will be used later on.  \end{remark}

 As an application of the above decomposition results,  we   give an  another proof of Proposition \ref{p:Jjump} (ii) and (iii).

\medskip
\noindent {\it Another proof of Proposition \ref{p:Jjump} (ii) and (iii)}. Notice that $T_{-m}^J=L(T_{-m}, 0)$. Conditioning on $T_{-m}^J=x$: by Corollary \ref{c:hitting2}, $x^{(m)}_1$ is distributed as $L(T_{J(x)}, 0)$ under $P^{\delta}$ biased by 
$J(x)^{ {\delta \over 2} -1 } {1}_{\{J(x)> -m \}}$. By the independence of $L(T_{J(x)}, 0)$ and $J(x)$ of Proposition \ref{p:min}, the biased law of $L(T_{J(x)}, 0)$ is the same as under $\p^\delta$, hence is $x$ times a beta(${\delta\over 2}, 1$) random variable. Moreover, conditionally on $x^{(m)}_1=y$ and $J(x)=-m_1$, the process before $T_{J(x)}$ is simply the process $X$ under $\p^\delta$ before hitting $T_{-m_1}$ conditioned on $L(T_{-m_1}, 0)=y$ (by Corollary \ref{c:hitting2} and Lemma \ref{l:perman}). Therefore we can iterate and get Proposition \ref{p:Jjump}. $\Box$

\section{Decomposition at the minimum} \label{s:decompositionminimum}
 
Let $\delta>0$.  Let $X_1$ be a copy of the process $X$ under $\p^{(-\delta)}$ and $X_2$  be an independent Bessel process of dimension 3, both starting at $0$.  Recall Notation \ref{n:<h}.  From our notations, $L^1(\infty, r)$, resp. $L^2(\infty, r)$, denotes the total local time at height $r$ of $X_1$, resp. $X_2$, while $X_1^{-, h}$, $X_2^{-,  h}$ are obtained by gluing the excursions   below $h$ of $X_1$ and $X_2$ respectively. We set
\begin{equation} \label{def:H12}
 H^{1,2} := \sup\{r\ge 0\, :\, L^1(\infty, r) + L^2(\infty, r) = 1 \}.
 \end{equation}

 Proposition \ref{p:legall-yor} yields that     the process $L^1(\infty, r) + L^2(\infty, r) , r\ge0$ is distributed as the square of a Bessel process of dimension $\delta+2$, starting from $0$. Then $H^{1,2} < \infty$ a.s.

 \begin{lemma}\label{l:XH} Let $\delta>0$. 
  Let $m>0$ and $x\in (0,1)$. Conditionally on $\{ H^{1,2}=m, L^1(\infty, H^{1,2}) =x \}$: 
\begin{enumerate}[(i)]
\item $X_1^{-,H^{1,2}}$ and $X_2^{-,H^{1,2}}$ are   independent; 
\item $X_1^{-,H^{1,2}}$ is distributed as $(X_t^{-,  m},\, t < A_\infty^{-, m})$ under $\p^{(-\delta)} ( \cdot | L^1(\infty, m) = x )$;
\item $X_2^{-,H^{1,2}}$ is distributed as $(X_t^{-,  m},\, t< A_\infty^{-, m})$ under $\p^{(-2)}( \cdot | L^2(\infty, m) = 1-x )$, 
\end{enumerate}
\noindent where  $ A_\infty^{-,  m}= \int_0^\infty 1_{\{ X_t \le m\}} {\rm d} t $ is the total lifetime of the process $X^{-,  m}$.
 \end{lemma}

  \medskip
  
 {\noindent\it Proof.}  First we describe the law of   $(X_t^{-,  m},\, t < A_\infty^{-, m})$ under $\p^{(-\delta)} ( \cdot | L^1(\infty, m) = x )$.   Let  ${\mathscr  D}_m:= \sup\{t>0: X_t \le m\}$ be the last passage time of $X$ at $m$ [note that under $\p^{(-\delta)}$, $X_t \to \infty$ as $t \to \infty$].  By the duality of Proposition \ref{p:duality},  $\{X_{{\mathscr  D}_m-t}- m, 0\le t \le {\mathscr  D}_m\}$, under $\p^{(-\delta)}$, has the same law as $\{X_t, 0\le t \le T_{-m}\}$ under  $\p^\delta$.     Corollary  \ref{c:hitting2}  gives  then  the law of $(X_t  , t \le {\mathscr  D}_m)$ under $\p^{(-\delta)} ( \cdot | L(\infty, m) = x )$.  The process $(X_t^{-, m},\, t < A_\infty^{-, m})$ is a measurable function of $(X_t  , t \le {\mathscr  D}_m)$. Note that $m \to (X_1^{-, m}, X_2^{-, m})$ is continuous. \footnote{For instance, we may show  that  for any   $T>0$, almost surely $\sup_{0\le t \le T}| X^{-,  m'}_t - X^{-, m}_t| \to 0$ as $m' \to m$. Let us give a proof  by contradiction. Suppose there exists some $\varepsilon_0>0$, a sequence $(t_k)$ in $[0, T]$ and $m_k \to m$ such that  $| X^{-, m_k}_{t_k} - X^{-,  m}_{t_k}|>\varepsilon_0$.  Write for simplification   $s_k:=\alpha_{t_k}^{-,  m}$ and $s'_k:=\alpha_{t_k}^{-,  m_k}$.  Consider the case $m_k >m$ (the other direction can be treated in a similar way). Then $s_k\ge s'_k$ and $|X_{s'_k}- X_{s_k}|> \varepsilon_0$  for all $k$.   Since $X_{s_k}\le m$ and $X_{s'_k}\le m_k$,   either $X_{s_k}\le m - \frac{\varepsilon_0}2$ or $X_{s'_k}\le m - \frac{\varepsilon_0}2$ for all large $k$.    Consider for example the case $X_{s'_k}\le m - \frac{\varepsilon_0}2$.  By the uniform continuity of $ X_t$  on every compact, there exists some $\delta_0>0$ such that $X_u\le m $ for all $|u-s'_k| \le \delta_0$ and $k\ge1$. Then for any $s \ge  s'_k$, $\int_{s'_k}^s  1_{\{X_u \le m\}} \d u    \ge \min(\delta_0, s-s'_k)$.  Note that by definition,   $\int_{s'_k}^{s_k} 1_{\{X_u \le m\}} \d u =t_k- \int_{0}^{s'_k} 1_{\{X_u \le m\}} \d u= \int_{m}^{m_k} L(s'_k, x) \d x \le \zeta\, (m_k- m)$, with $\zeta:=\sup_{x\in \r} L(\alpha_T^{-,m}, x)$.   It follows that for all sufficiently large $k$, $0\le s_k-s'_k \le \zeta \, (m_k- m)$. Consequently $X_{s_k}- X_{s'_k} \to 0$ as $m_k \to m$, in contradiction with the assumption that $|X_{s'_k}- X_{s_k}|> \varepsilon_0$  for all $k$. This proves the continuity of $m \to X^{-,  m}$.  }   From Corollary \ref{c:hitting2}, we may find a regular version of the    law  of   $(X_t^{-, m},\, t < A_\infty^{-, m})$ under $\p^{(-\delta)} ( \cdot | L(\infty, m) = x ) $ such that    for any bounded continuous functional  $F$ on  ${\mathcal K}$,  the application $$(m, x) \mapsto  \e^{(-\delta)} \big[ F(X_t^{-, m},\, t < A_\infty^{-, m}) | L(\infty, m) = x  \big]$$ is continuous.

  Now let us write $H:= H^{1, 2}$ for concision. 
  Let $F_1$ and $F_2$ be two bounded  continuous functionals on    ${\mathcal K}$  and $g: \r^2_+ \to \r$ be a bounded continuous function.   Let $H_n:= 2^{-n} \lfloor 2^n H \rfloor$ for any $n\ge1$.  By the continuity of $m \to (X_1^{-, m}, X_2^{-, m})$ and  that of $(L^1(\infty, m), L^2(\infty, m))$, we have $$ \e\Big [ F_1(X_1^{-, H}) F_2(X_2^{-, H}) g(H, L^1(\infty, H))\Big]
  =
  \lim_{n \to \infty} \e\Big [ F_1(X^{1, H_n}) F_2(X^{2, H_n}) g(H_n, L^1(\infty, H_n))\Big].$$

  Note that \begin{eqnarray*}
  &&
   \e\Big [ F_1(X_1^{-, H_n}) F_2(X_2^{-, H_n}) g(H_n, L^1(\infty, H_n))\Big]
  \\
   &=&
   \sum_{j=0}^\infty \e\Big [ F_1(X_1^{-, \frac{j}{2^n}}) F_2(X_2^{-, \frac{j}{2^n}}) g\big(\frac{j}{2^n}, L^1(\infty, \frac{j}{2^n})\big) 1_{\{ \frac{j}{2^n} \le H < \frac{j+1}{2^n}\}} \Big].
   \end{eqnarray*}

 By the independence property of Corollary \ref{c:ind} and the duality of Proposition \ref{p:duality}, conditioning on  $\{L^1(\infty, \frac{j}{2^n}), L^2(\infty, \frac{j}{2^n})\}$, the processes $(X_1^{-, \frac{j}{2^n}}, X_2^{-, \frac{j}{2^n}})$ are independent,  and independent of $(X_1^{+,   \frac{j}{2^n}}, X_2^{+,   \frac{j}{2^n}})$. Since $\{\frac{j}{2^n} \le H < \frac{j+1}{2^n}\}$ is measurable with respect to $\sigma(X_1^{+,  \frac{j}{2^n}}, X_2^{+,    \frac{j}{2^n}})$, we get that for each $j\ge0$, \begin{eqnarray*}
  && \e\Big [ F_1(X_1^{-, \frac{j}{2^n}}) F_2(X_2^{-, \frac{j}{2^n}}) g\big(\frac{j}{2^n}, L^1(\infty, \frac{j}{2^n})\big) 1_{\{ \frac{j}{2^n} \le H < \frac{j+1}{2^n}\}} \Big]
  \\
  &=&
  \e \Big[\Phi_1\big(\frac{j}{2^n}, L^1(\infty, \frac{j}{2^n})\big) \, \Phi_2\big(\frac{j}{2^n}, L^2(\infty, \frac{j}{2^n})\big) \,  g\big(\frac{j}{2^n}, L^1(\infty, \frac{j}{2^n})\big) 1_{\{ \frac{j}{2^n} \le H < \frac{j+1}{2^n}\}} \Big],
   \end{eqnarray*}

   \noindent  where $$ \Phi_1(m, x):= \e[ F_1(X_1^{-,  m})\, | L^1(\infty, m)=x], \, \,  \Phi_2(m, x):= \e[ F_2(X_2^{-,  m})\, | L^2(\infty, m)=x].$$

 \noindent By Remark \ref{r:Phigm}, $\Phi_1$ and $\Phi_2$ are continuous functions in $(m, x)$.   Taking the sum over  $j$ we get that \begin{eqnarray*}
  &&
   \e\Big [ F_1(X_1^{-, H_n}) F_2(X_2^{-, H_n}) g(H_n, L^1(\infty, H_n))\Big]
  \\
   &=&
\e\Big [ \Phi_1( H_n,  L^1(\infty, H_n)) \Phi_2( H_n, L^2(\infty, H_n)) g(H_n, L^1(\infty, H_n))\Big] .
   \end{eqnarray*}

 \noindent Since $\Phi_1$ and $\Phi_2$ are bounded and continuous,  the dominated convergence theorem yields  that \begin{eqnarray}
&&  \e\Big [ F_1(X_1^{-, H}) F_2(X_2^{-, H}) g(H, L^1(\infty, H))\Big]
\nonumber  \\
 & =&
  \lim_{n \to \infty} \e\Big [ F_1(X_1^{-, H_n}) F_2(X_2^{-, H_n}) g(H_n, L^1(\infty, H_n))\Big]
 \nonumber  \\
  &=&
   \e\Big [ \Phi_1( H ,  L^1(\infty, H)) \Phi_2( H , L^2(\infty, H)) g(H , L^1(\infty, H))\Big] ,  \label{F12PHI12}
 \end{eqnarray}
 
 \noindent  proving  Lemma \ref{l:XH} as $L^2(\infty, H) = 1- L^1(\infty, H).$ $\Box$

\medskip

\begin{remark} \label{r:XH} Let $\delta>2$. Consider the process $X$ under $\p^{(-\delta)}$ and the total  local time  $L(\infty, r)$ of $X$ at position $r\ge0$.  For $x>0$, let   $$H_x:= \sup\{r \ge0: L(\infty, r) =x\}.$$

\noindent By  Proposition \ref{p:legall-yor},  $H_x < \infty$, $\p^{(-\delta)}$-a.s. Note  that the same arguments  leading to \eqref{F12PHI12}     shows that  $$ \e\Big [ F_1(X^{-,  H_x})   g(H_x, L(\infty, H_x))\Big] =  \e\Big [ \Phi_1( H_x ,  L(\infty, H_x))   g(H_x , L(\infty, H_x))\Big], $$

\noindent where $\Phi_1(m, x):= \e[ F_1(X^{-,  m})\, | L(\infty, m)=x]$ for $m>0, x> 0$.  Since $L(\infty, H_x)=x$, we obtain that   conditionally on $\{H_x= m\}$,  the process $X^{-,  H_x}$ is distributed as $(X_t^{-,  m},\, t < A_\infty^{-,  m})$ under $\p^{(-\delta)} ( \cdot | L(\infty, m) = x )$.
\end{remark}

 Recall  \eqref{def-T}, \eqref{def-Jx} and \eqref{def:H12}. The main result in this section is the following path decomposition of $(X_t)$ at $T_{J(1)} = \inf\{t \in [0, \tau_0(1)] : X_t= J(1)\}$, the unique time before $\tau_0(1)$ at which $X$ reaches  its minimum $J(1)$.
 
 \begin{theorem}  \label{p:decomposition} Let $\delta>0$.   Define $Z_1:= (X_{T_{J(1)} - t} - J(1))_{t\in [0,T_{J(1)} ]}$ and $Z_2:=  (X_{T_{J(1)} + t} - J(1))_{t\in [0,\tau_0(1) - T_{J(1)} ]}$. Under $\p^\delta$, the couple of processes 
$$
 \Big( Z_1^{-,  |J(1)|}, \, Z_2^{-,  |J(1)|} \Big)$$

\noindent    is distributed  as  $(X_1^{-,H^{1,2}},X_2^{-,H^{1,2}})$.

\end{theorem}

\noindent {\it Proof}.  From Lemmas \ref{l:perman} and \ref{l:XH}, it remains to prove that the joint law of $(|J(1)|,L(T_{J(1)}, 0))$ is the same as $(H^{1,2},L^1(\infty, H^{1,2}))$. Recall the law of $(J(1) , L(T_{J(1)}, 0))$ from Propositions  \ref{p:Jjump} (i) and \ref{p:min}. Define a process  $(Y_t,\, -\infty < t <\infty)$ with values in $[0,1]$ defined by time-change as 
$$
Y_{ A_m} = {L^1(\infty, m) \over L^1(\infty, m) + L^2(\infty, m)} , 
$$
 
\noindent where $A_m := \int_1^m {{\rm d}h \over L^1(\infty, h) + L^2(\infty, h)}$ for any $m>0$ (as such $  \lim_{m\to 0} A_m = - \infty$ a.s.). Following Warren and Yor \cite{warren-yor}, equation (3.1), we call Jacobi process of parameters $d,d' \ge 0$ the diffusion with generator $2y(1-y) \frac{\d^2}{\d y^2} + (d - (d+d')y) \frac{\d}{\d y}$. We claim that $Y$ is a stationary Jacobi process of parameter $(\delta , 2)$, independent of $(L^1(\infty, m)+L^2(\infty, m),\, m\ge 0)$.  It is a consequence of Proposition 8 of Warren and Yor \cite{warren-yor}.  Let us see why.

First, notice that $  {L^1(\infty, m) \over L^1(\infty, m) + L^2(\infty, m)}$ is a beta$({\delta \over 2},1)$-random variable for any $m>0$, because $ L^1(\infty, m) $ and $L^2(\infty, m)$ are independent and distributed as gamma(${\delta \over 2},2m$) and  gamma$(1,2m)$ respectively, by Proposition \ref{p:legall-yor} and the duality in Proposition \ref{p:duality}.  It is independent of $L^1(\infty, m)+L^2(\infty, m)$, hence by the Markov property, also of $(L^1(\infty, h)+L^2(\infty, h),\, h\ge m )$.

Let $t_0\in \r$. By Proposition 8 of   \cite{warren-yor}, for any $m\in (0, 1)$, conditioning on $(L^1(\infty, h)+L^2(\infty, h),\, h\ge m)$ , the process $(Y_{h+A_m},\, h\ge 0)$ is distributed as a Jacobi process starting from a beta$(\frac\delta2, 1)$ random variable, hence stationary. Notice that $A_m$ is measurable with respect to $\sigma( { L^1(\infty,  h) + L^2(\infty, h)}, \, h\ge m)$. We deduce that, conditioned on  $(L^1(\infty, h)+L^2(\infty, h),\, h\ge m)$ and $A_m\le t_0$, the process $(Y_h,\, h \ge t_0)$ is a Jacobi process starting from a beta(${\delta \over 2},1$) random variable. Letting  $m\to 0$ we see that $Y$ is a stationary Jacobi process of parameter $(\delta , 2)$, independent of $(L^1(\infty, m)+L^2(\infty, m),\, m\ge 0)$.

Since $L^1(\infty, H^{1,2})= Y_{A_{H^{1,2}}}$ and $A_{H^{1,2}}, H^{1,2}$ are measurable with respect to $\sigma\{L^1(\infty, m)+L^2(\infty, m),\, m\ge 0\}$, we  deduce that the random variable $L^1(\infty, H^{1,2})$ follows the beta(${\delta\over 2},1$) distribution and that $H^{1,2}$ and $L^1(\infty,H^{1,2})$ are independent.   Finally, the random variable $H^{1,2}$ is the exit time  at $1$ of a square Bessel process of dimension $2+\delta$ by Proposition \ref{p:legall-yor}, whose density is equal to $\frac1{\Gamma(\frac\delta2)} 2^{-\frac\delta2} t^{-\frac\delta2-1} e^{-\frac1{2t}}$ for $t>0$ (Exercise (1.18), Chapter XI of Revuz and Yor \cite{revuz-yor}). By  Proposition \ref{p:Jjump} (i), we see that $|J(1)|$   is distributed as  $H^{1,2}$. This completes the proof.  $\Box$

\medskip
The rest of this section is devoted to  a path decomposition of $X$ under $\p^{(-\delta)}$ for $\delta>2$. 
For $x>0$, let as in Remark \ref{r:XH},
$$
H_x :=\sup\{  r \ge 0 \,:\, L(\infty, r) =x\}.  
$$

\noindent Define  $$
  S_x:=\sup\{ t\ge 0\,:\, X_t=H_x\} , \qquad  \widehat J_x:=\inf\{X_t,\, t\ge T_{H_x}\}-H_x , 
$$
where as before  $T_{H_x}:= \inf\{t\ge0: X_t= H_x\}$ is the hitting time of $H_x$ by $X$.

Write $C_x:=H_x+ \widehat J_x$. We consider the following three processes: 
\begin{eqnarray*}
 X^{(1)} &:=&   (X_{S_x -t} - H_x,\, t\in [0,S_x-T_{H_x}]), \\
 X^{(2)} &:=& - (X_{T_{H_x}-t}-H_x,\, t \in [0,T_{H_x}- D_x]), \\
 X^{(3)} &:=& (X_{D_x-t}-C_x,\, t \in [0, D_x]),
\end{eqnarray*}

\noindent where $D_x:= \sup\{t< T_{H_x}: X_t= C_x\}$.

Furthermore,   let   $X^{(1), -}$ be the process $X^{(1)}$ obtained by  removing all its  positive excursions: $$X^{(1), -}_t:= X^{(1)}_{\alpha_t^{(1), -}},  $$ with $\alpha_t^{(1), -}:= \inf\{s>0: \int_0^s 1_{\{ X^{(1)}_u  \le 0\}}  \d u$  and $t \le \int_0^{S_x-T_{H_x}} 1_{\{ X^{(1)}_u  \le 0\}}  \d u.$

\begin{proposition} \label{p:mucond} 
Let $ \delta > 2 $ and $ x ,  a>0$.  

(i)  Under $\p^{(-\delta)}$, $\frac1{|\widehat J_x|}$ is distributed as ${\mbox gamma}(\frac\delta2, \frac{2}{x})$. 

(ii) Under $\p^{(-\delta)} (\cdot \, |\, \widehat J_x= -a)$,  the three processes $X^{(1), -}$, $X^{(2)}$, $X^{(3)}$ are independent and distributed respectively as 

$\bullet$   $(X_t, 0\le t \le  \tau_0(x))$, under $\p^\delta( \cdot |   J(x) =-a)  $,  after removing all excursions above $0$;

$\bullet$   a Bessel process $(R_t)_{0\le t \le T_a} $of dimension 3 starting from $0$ killed at $T_a:= \inf\{t>0: R_t=a\}$;

$\bullet$    $(X_t, 0\le t \le T_{-a(1-u)/u})$  under $\p^\delta (\cdot | T_{-a(1-u)/u} < T_a)$, where $u \in [0, 1]$ is independently chosen  according to the law beta(${\delta\over 2}-1,1$).
 
\end{proposition}

Since $\widehat J_x$ under $\p^{(-\delta)}$ is distributed as $J(x)$ under $\p^\delta$, we observe that the (unconditional) law of $X^{(1), -}$ under $\p^{(-\delta)} $ is equal to that    $(X_t, 0\le t \le  \tau_0(x))$, under $\p^\delta$,  after removing all excursions above $0$.  

\medskip
{\noindent\it Proof.} Let $m>0$.   By Remark \ref{r:XH},  conditionally on $\{H_x= m\}$,  the process $X^{-,  H_x}$ is distributed as $(X_t^{-,  m},\, t < A_\infty^{-,  m})$ under $\p^{(-\delta)} ( \cdot | L(\infty, m) = x )$.   

Note that  conditionally on $\{H_x= m\}$,    $S_x= \sup\{t>0: X_t \le m\}$ is  the last passage time of $X$ at $m$.   By the duality of Proposition \ref{p:duality},  under $\p^{(-\delta)} ( \cdot | L(\infty, m) = x )$,  the process $(X_{A_\infty^{-,  m}-t}^{-,  m}-m,\, t < A_\infty^{-,  m})$ is distributed as    $\{X^{-,  0}_t, 0\le t \le A_{T_{-m}}^{-,  0}\}$ under  $\p^\delta( \cdot | L(T_{-m}, 0)=x)$, the process $(X_t, 0\le t \le T_{-m})$ obtained by  removing all positive excursions.     Furthermore,  remark that $\widehat J(x)$ corresponds to $I_{{\tt g}_m}$ which is defined for  the process $(X_t, 0\le t \le T_{-m})$ under $\p^\delta( \cdot | L(T_{-m}, 0)=x)$.   Then $\p^{(-\delta)}( | \widehat J_x| \in \cdot   \, |\, H_x=m) = \p^\delta( |I_{{\tt g}_m} |\in \cdot \, |\, L(T_{-m}, 0)=x)$. We deduce that for any $0< a <m$, the conditional law of the process $(X_{A_\infty^{-,  m}-t}^{-,  m}-m,\, t < A_\infty^{-,  m})$  under $\p^{(-\delta)}( \cdot \,|\,  H_x=m, \widehat J_x=-a)$ is the same as the conditional law of the process $\{X^{-,  0}_t, 0\le t \le A_{T_{-m}}^{-,  0}\}$ under  $\p^\delta( \cdot | L(T_{-m}, 0)=x, I_{{\tt g}_m} =-a).$

Then we may apply  Theorem \ref{p:hitting2}  (ii)   to see that   conditionally on $\{H_x= m, \widehat J_x= -a \}$, $X^{(1), -}$, $X^{(2)}$, and  $X^{(3)}$ are independent, and 

$\bullet$ $X^{(1), -}$ is distributed as $(X_t^{-,  0}, t \le A^{-,  0}_{\tau_0(x)})$  under $\p^\delta( \cdot |   J(x) =-a)  $,  where $(X_t^{-,  0}, t \le A^{-,  0}_{\tau_0(x)})$ is the process obtained from   $(X_t, 0\le t \le  \tau_0(x))$ by removing all positive excursions;

$\bullet$ $X^{(2)}$ is distributed as a three-dimensional Bessel process $(R_t)_{0\le t \le T_a}$ killed at $T_a:= \inf\{t>0: R_t=a\}$;

$\bullet$ $X^{(3)}$ is distributed as  $(X_t, 0\le t \le T_{-(m-a)})$  under $\p^\delta (\cdot | T_{-(m-a)} < T_a)$. 

   Moreover   \begin{eqnarray*}  \p^{(-\delta)} (|\widehat J_x| \in \d a \, |\, H_x  =  m) 
 &=& \p^\delta(I_{{\tt g}_m} \in \d a \, |\, L(T_{-m}, 0)=x)
 \\
 &=&  \frac{x}{2 a^2} 
e^{-\frac{x}{2a}  + \frac{x}{2 m}} 1_{\{0< a < m\}}  \d a ,
\end{eqnarray*}

\noindent  where the last equality follows from   \eqref{l:Jcond}.   

Recall (\cite{GY03}) the law of $H_x$ under $\p^{(-\delta)}$ :  For $\delta>2$,  $$ \p(H_x \in \d m) / \d m = \frac{(\frac{x}2)^{\frac\delta2-1}}{\Gamma(\frac\delta2-1)} m^{-\frac\delta2} \, e^{-\frac{x}{2m}}, \qquad m>0.$$

\noindent We get  $$\p^{(-\delta)} (|\widehat J_x| \in \d a)
=
\frac{1}{\Gamma(\frac\delta2)} (\frac{x}{2})^{\frac\delta2} \, a^{-(\frac\delta2+1)}\, e^{-\frac{x}{2 a}} \, \d a, \qquad a>0,$$

\noindent which implies (i).

For any bounded continuous functionals $F_1, F_2, F_3$ on ${\mathcal K}$, we have    \begin{eqnarray*}
&& \e^{(-\delta)} [ F_1(X^{(1, \le 0)})\, F_2(X^{(2)})\, F_3 (X^{(3)})\, |\, \widehat J_x= -a]
  \\
&=&
\int _a^\infty  \frac{\p^{(-\delta)}(|\widehat J_x| \in \d a, H_x \in \d m)}{\p(|\widehat J_x| \in \d a)} \e^{(-\delta)} [ F_1(X^{(1), -}\, F_2(X^{(2)})\, F_3 (X^{(3)})\, |\, \widehat J_x= -a, H_x=m ] 
  \\
&=&
(\frac\delta2-1)\,  \int_a^\infty \d m  a^{\frac\delta2-1} m^{-\frac\delta2}  \,  \e^\delta [ F_1(X_t^{-,  0}, t \le A^{-,  0}_{\tau_0(x)}) | J(x)=-a]\, \e [F_2( R_t, t\le T_a)]\, 
  \\
&& \qquad  \times \e^\delta [  F_3 ( X_t,  t \le T_{-(m-a)})\, |\,T_{-(m-a)} < T_a ]   
\\
&=&   \e^\delta [ F_1(X_t^{-,  0}, t \le A^{-,  0}_{\tau_0(x)}) | J(x)=-a]\, \e [F_2( R_t, t\le T_a)]\,  \times
\\
&& \quad (\frac\delta2-1)\, \int_0^1 \d u u^{\frac\delta2-2}   \,    \e^\delta [  F_3 ( X_t,  t \le T_{-a(1-u)/u})\, |\,T_{-a(1-u)/u} < T_a ] ,
\end{eqnarray*}

\noindent which gives  (ii) and completes the proof of the Proposition.  $\Box$

\section{The perturbed Bessel process and  its rescaling at a stopping time}\label{s:bessel}

We rely on the paper of Doney, Warren and Yor \cite{dwy98},  restricting  our attention to the case of dimension $d=3$. For $\kappa<1$, the $\kappa$-perturbed Bessel process of dimension $d=3$ starting from $a\ge 0$ is the process $(R_{3, \kappa},\, t\ge 0)$ solution of 
\begin{equation}\label{R3alpha}
R_{3, \kappa}(t)  = a + W_t +  \int_0^t {  {\rm d} s\over R_{3, \kappa}(s) } + \kappa ( S_t^{R_{3, \kappa}} - a),  
\end{equation}

\noindent where $S_t^{R_{3, \kappa}} = \sup_{0\le s\le t} R_{3, \kappa}(s)$ and $W$ is a standard  Brownian motion. For $a>0$, it can be constructed as the law of $X$ under the measure $\p^{3,\kappa}_a$ defined by
\begin{equation}\label{eq:radon}
\p^{3,\kappa}_a\big|_{{\cal F}_t} = {1\over a^{1-\kappa}} { X_{t\land T_0} \over (S_{t\land T_0})^\kappa } \tilde \p^\delta_a\big|_{{\cal F}_t}
\end{equation}

\noindent where: $\delta:=2(1-\kappa)$,   $S_t:=\sup_{0\le s\le t} X_s$, and for any $a\ge 0$, $X$ under $\tilde \p^\delta_a$ is distributed as $a-X$ under $\p^\delta$.   Roughly speaking, the $\kappa$-perturbed Bessel process of dimension $3$ can be thought of as the process $-X$ under $\p^{\delta}$ ``conditioned to stay positive''. The next proposition is very related to Lemma 5.1 of \cite{dwy98}.  We let $\p^{3,\kappa}=\p^{3,\kappa}_0$ be a probability measure under which $X$ is distributed as the process $R_{3,\kappa}$ starting from $0$. Under $\p^{3,\kappa}_a$ and $\tilde \p^{\delta}_a$,  we denote by $\{ (q, {\mathfrak e}_{X,q}^\downarrow),\, q \in {\mathcal Q}^\downarrow_X\}$   the excursions of the process   $ (S_t,X_t-S_t )$ away from  $\r \times \{0\}$. Notice that, under $\tilde \p^{\delta}_0$, by Proposition \ref{p:loopsoup} and the invariance in distribution of the loop soup by the map $x\to -x$,  these excursions are distributed as $\{ (q, {\mathfrak e}_q^\downarrow),\, q \in {\mathcal Q}^\downarrow_{\delta \over 2} \cap (0,\infty)\}$.

\begin{proposition}\label{p:besselloop}
Let $\kappa<1$ and $\delta:=2(1-\kappa)$. The point process $\{(q,{\mathfrak e}_{X,q}^\downarrow),\, q \in {\mathcal Q}_X^\downarrow\}$ under $\p^{3,\kappa}$ is distributed as the Poisson point process 
$$
\{(q,{\mathfrak e}_q^\downarrow),\, q \in {\mathcal Q}^{\downarrow}_{\delta \over 2} \hbox{  such that }    q+ {\mathfrak e}_q^\downarrow  \subset (0, \infty)
\}.
$$ 

\noindent   In other words, the excursions of $R_{3, \kappa}$ below its supremum, seen as unrooted loops, are distributed as the loops of the Brownian loop soup ${\mathcal L}_{\delta\over 2}$ which entirely lie in the positive half-line.
\end{proposition}

\begin{remark} \label{r:R3alpha}   

(i) The intensity measure of this Poisson point process has been computed in \eqref{down} of Lemma \ref{l:intensity}.

(ii)  Similarly to Section 5.1 of Lupu \cite{lupu}, one can construct the process $R_{3, \kappa}$ from the loops of the Brownian loop soup ${\mathcal L}_{\delta\over 2}$ which entirely lie in $(0,\infty)$, by rooting them at their maxima and gluing them in the increasing order of their maxima.   

 (iii)  The process  $(R_{3, \kappa}, S^{R_{3, \kappa}})$ is a Markov process. Hence, applying the strong Markov property under $\p^{3,\kappa}$ to $X$ at time $T_a$, we deduce that under $\p^{3,\kappa}_a$, the excursions below supremum of the process, seen as unrooted loops, are distributed as the loops of the Brownian loop soup ${\mathcal L}_{\delta\over 2}$ which entirely lie in the positive half-line and with maximum   larger than  $a$. It entails that for $a>0$, the process $X$ under $\p^{3,\kappa}_a$ is the limiting distribution as $m\to \infty$ of the process $X$ under $\tilde \p^\delta_a$ conditioned on hitting $m$ before $0$ (the process $X$ is measurable with respect to its excursions below supremum, which are equally distributed before time $T_m$ under $\p^{3,\kappa}_a$ and under $\tilde \p^\delta_a( \cdot | T_m <T_0)$).

\end{remark}

\noindent {\it Proof of Proposition \ref{p:besselloop}}.    Let $f: \r_+\times {\mathcal K} \to \r_+$  be    measurable. For any $0< s<s'$, we compute 
 $$
\e^{3, \kappa} \Big[\ee^{- \sum_{q\in {\mathcal Q}_X^\downarrow \cap [s,s'] }  f(q, {\mathfrak e}_{X,q}^\downarrow)} \Big].
$$
\noindent Notice that the integrand is measurable with respect to the $\sigma$-algebra   $\sigma(X_t,\, t\in [T_s,T_{s'}])$.  By the strong Markov property at time  $T_s$ and the absolute continuity \eqref{eq:radon} with $a=s$ there, the previous expectation is equal to
$$
s^{\kappa-1} \frac{s'}{(s')^\kappa} \, \widetilde \e^{\delta}_0 \Big[ \ee^{- \sum_{q\in {\mathcal Q}_X^\downarrow \cap [s,s'] }  f(q, {\mathfrak e}_{X,q}^\downarrow)} \,  1_{\{ T_0 \circ \theta_{T_s} > T_{s'}\}} \Big]
$$

\noindent where $\theta$ is the shift operator.  Notice that  
$$
  \ee^{- \sum_{q\in {\mathcal Q}_X^\downarrow \cap [s,s'] }  f(q, {\mathfrak e}_{X,q}^\downarrow)} \,  1_{\{ T_0 \circ \theta_{T_s} > T_{s'}\}}
=  \ee^{- \sum_{q\in {\mathcal Q}_X^\downarrow \cap [s,s'] , q+{\mathfrak e}_{X,q}^\downarrow  \subset (0, \infty)}  f(q, {\mathfrak e}_{X,q}^\downarrow)} \, 1_{\mathcal{E}}
$$

\noindent where $\mathcal E$ is the event that the set of  $q\in {\mathcal Q}_X^\downarrow \cap [s,s']$ such that $q+{\mathfrak e}_{X,q}^\downarrow  \not\subset (0, \infty)$ is empty.   
  
 We already mentioned that  the collection of $(q,{\mathfrak e}_{q,X}^\downarrow)$ for $q\in {\mathcal Q}_X^\downarrow$ is a Poisson point process under $\tilde \p^{\delta}_0$, distributed as $\{ (q, {\mathfrak e}_q^\downarrow),\, q \in {\mathcal Q}^\downarrow_{\delta \over 2} \cap (0,\infty)\}$. By the independence property of Poisson point processes, we deduce that 
$$
\e^{3, \kappa} \Big[\ee^{- \sum_{q\in {\mathcal Q}_X^\downarrow \cap [s,s'] }  f(q, {\mathfrak e}_{X,q}^\downarrow)} \Big]
=
c   \widetilde \e^{\delta} \Big[  \ee^{- \sum_{q\in {\mathcal Q}_X^\downarrow \cap [s,s'] , q+{\mathfrak e}_{X,q}^\downarrow  \subset (0, \infty)}  f(q, {\mathfrak e}_{X,q}^\downarrow)}      \Big]
$$

\noindent for some constant $c$ which is necessarily $1$. The Proposition follows.  $\Box$

\medskip
In the rest of this section,  we will extend  the definition of perturbed Bessel processes to allow some positive local time at $0$.

Let $x\ge 0$.  
We define a kind of {\it perturbed Bessel process  $R_{3, \kappa}^x$  with local time $x$ at position $0$}.  More precisely, for $\kappa<1$ and $\delta=2(1-\kappa)$, we denote by $R_{3, \kappa}^x$ the process obtained  by concatenation in the following way: take $-X$ under $\p^{\delta}$ up to time $\tau_0(x)$, biased by $|J(x)|^{{\delta\over 2}-1}$ followed by  a Bessel of dimension 3 killed when hitting $|J(x)|$, followed by the  $\kappa$-pertubed process $R_{3, \kappa}$ starting from $|J(x)|$.  Recall Theorem \ref{p:hitting2}. Corollary \ref{c:hitting2} and Remark \ref{r:R3alpha} (iii)  show that, when $x>0$, $R_{3, \kappa}^x$ is the limit in distribution of the process $(-X_t,\, t\le T_{-m})$ under $\p^\delta(\cdot | L(T_{-m}, 0)=x)$ as $m\to \infty$.   Clearly when $x=0$, $R_{3, \kappa}^0$   coincides with $R_{3, \kappa}$ defined previously in \eqref{R3alpha}  with $a=0$.  The following theorem shows that one can recover the process $R_{3,\kappa}^x$ by a suitable time-space scaling of a conditioned PRBM up to a hitting time. It is an extension of Theorem 2.2 equation (2.5) of Doney, Warren, Yor  \cite{dwy98} [which  corresponds  to the case $x=0$ and $m=1$].

\begin{theorem}\label{p:scaling1'}
Suppose $\kappa<1$ and let $\delta:=2(1-\kappa)$. Fix $m > 0$ and $x\ge 0$.   Let the space-change
$$
\theta(z) := 
\left\{
\begin{array}{ cc}
- {m z \over  m+z} & {\rm if} \; z\ge 0, \\
-z & {\rm if} \; z<0, 
\end{array} 
\right.  
$$
and the time-change
$$
A_t := \int_0^t \big(\theta'\big(R^x_{3,\kappa}(s)\big)\big)^2 \d s, \qquad t\ge 0.
$$

\noindent If $\widetilde X$ is defined via $\theta\big(R^x_{3,\kappa}(t)\big):= \widetilde X_{A_t} $, then $\widetilde X$ is distributed as $(X_t, 0\le t \le T_{-m})$ under $\p^\delta(\cdot | L(T_{-m}, 0)=x)$.
\end{theorem}

{\noindent\it Proof}.  
First we describe the  excursions above infimum of $X$  under $\p^\delta(\cdot | L(T_{-m}, 0)=x)$ in terms of the Brownian loop soup ${\mathcal L}_{\delta\over 2}$.  For the loops which hit $0$, we use again the same observation  as   in the proof of  \eqref{l:Jcond}:    
conditionally on $\{\ell_\gamma^0:  \gamma \in {\mathcal L}_{\delta\over 2}, 0 \in \gamma\} $, the loops $\gamma$ such that   $\min\gamma>-m$ are (the projection on the space of unrooted loops of)  independent Brownian motions stopped at local time given by $\ell_\gamma^0$, conditioned on not hitting $-m$.  

Remark that  the set $\{\ell_\gamma^0:  \gamma \in {\mathcal L}_{\delta\over 2}, 0 \in \gamma, \min \gamma >- m\} $   is equal to the (non-ordered) set $\{x_{i-1}^{(m)}- x_i^{(m)}\}_{i\ge1}$ in the notation of Proposition \ref{p:Jjump},    and $L(T_{-m}, 0)= T^J_{-m}$.  By Proposition \ref{p:Jjump},  conditionally  on $\{L(T_{-m}, 0)=x\}$,  the ordered sequence of  $(\ell^0_\gamma: \gamma \in  {\mathcal L}_{\delta\over 2} , 0\in \gamma, \min \gamma > -m)$ is distributed as  $x(P_{(1)},P_{(2)},\ldots)$, where $(P_{(1)},\, P_{(2)}, \ldots)$ has the Poisson--Dirichlet distribution of parameter ${\delta\over 2}$.

Note that the loops $\gamma$ of ${\mathcal L}_{\delta\over 2}$  such that $  \gamma \subset  (-m, 0)$ are independent of $L(T_{-m}, 0)$. 
Then  the  excursions above infimum of $X$ (seen as unrooted loops)  under $\p^\delta(\cdot | L(T_{-m}, 0)=x)$ consists in the superposition of:
\begin{itemize}
\item the loops $\gamma$ of the Brownian loop soup ${\mathcal L}_{\delta\over 2}$  such that $  \gamma \subset  (-m, 0)$;
\item an independent collection of independent Brownian motions stopped at local time $x(P_{(1)},P_{(2)},\ldots)$, where $(P_{(1)},\, P_{(2)}, \ldots)$ has the Poisson--Dirichlet distribution of parameter ${\delta\over 2}$.
\end{itemize}

Let $m\to \infty$,  we deduce that the excursions below supremum of $R_{3, \kappa}^x$ (seen as unrooted loops) consists in the superposition of:
\begin{itemize}
\item the loops $\gamma$ of the Brownian loop soup ${\mathcal L}_{\delta\over 2}$  such that $\min \gamma>0$;
\item an independent collection of independent Brownian motions stopped at local time $x(P_{(1)},P_{(2)},\ldots)$, where $(P_{(1)},\, P_{(2)}, \ldots)$ has the Poisson--Dirichlet distribution of parameter ${\delta\over 2}$.
\end{itemize}

Note  that via the transformation $\theta(R^x_{3,\kappa}(t))=\widetilde X_{A_t}$, the excursions of $R^x_{3,\kappa}$ below their current supremum are transformed into excursions of $\widetilde X$ above their current infimum. Since the law of $R^x_{3,\kappa}$ is characterized by the law of its excursions below their current supremum (exactly as Remark \ref{r:R3alpha} (ii)),  and the law of $\widetilde X$ is characterized by the law of its excursions above  their  current infimum in a similar way, we only have to focus on the law of those excursions and show that applying the transformation in space $\theta$ and in time $A_t^{-1}$, say $\Phi$,\footnote{More precisely for any process $(\gamma_t, t\ge0)$, $\Phi(\gamma)$ is the process defined by $\theta(\gamma_t)= \Phi(\gamma)\big(\int_0^t (\theta'(\gamma_s))^2 \d s\big)$. }

(a)  the loops $\gamma$ of the Brownian loop soup ${\mathcal L}_{\delta\over 2}$  such that $\min \gamma>0$ are transformed  into loops $\widetilde \gamma $ of ${\mathcal L}_{\delta/2}$ such that $\max \widetilde \gamma <0$ and $\min \widetilde \gamma > -m$;

(b)  for any $\ell>0$,  a Brownian motion $(B_t, 0\le t \le \tau_\ell^B)$ stopped at local time $\ell$ is transformed   into $(B_t, 0\le t \le \tau_\ell^B)$ conditioned on $\{\inf_{0\le t\le  \tau_\ell^B} B_t > -m\}$.

Let us prove (a).  For a loop $\gamma$, we let $\gamma^{\uparrow}$ be the loop $\gamma-\min\gamma$ rooted at its minimum,  and $\gamma^{\downarrow}$ be the loop $\gamma-\max\gamma$ rooted at its maximum. Notice that $\gamma^\uparrow$ is a positive excursion above $0$, and $\gamma^{\downarrow}$ is a negative excursion below $0$. We remark that for any loop $\gamma$ with $\min \gamma >0$,   $  \min \Phi(\gamma) = \theta(a) $ with $a:=\max \gamma$,    and  $\Phi(\gamma)^\uparrow= \Phi(a+ \gamma^\downarrow)- \theta(a)$.  By Lemma \ref{l:intensity}    \eqref{down},   for any nonnegative measurable function $f$ on $\r_- \times {\mathcal K}$, we have \begin{eqnarray*}
&&
\e \Big[ e^{-\sum_{  \gamma\in{\mathcal L}_{\delta/2}, \,  \min \gamma>0} f(  \min \Phi(\gamma)  ,   \Phi(\gamma)^\uparrow)}\Big]
\\
&=&
\exp\Big(-\delta  \int_0^\infty \d a \int    {\mathfrak n}^- (\d \mathfrak e) (1- e^{-f(\theta(a), \Phi(a +{\mathfrak e})-\theta(a))}) 1_{\{  \min {\mathfrak e} >-  a\}}\Big)\\
&=&
\exp\Big(-\delta  \int_0^\infty \d a \int    {\mathfrak n}^+ (\d \mathfrak e) (1- e^{-f(\theta(a), \Phi(a-{\mathfrak e})-\theta(a))}) 1_{\{  \max{\mathfrak e} <  a}\}\Big)
.
\end{eqnarray*}

Let $h>0$.     Williams' description of the It\^{o} measure says that under ${\mathfrak n}^+ (\cdot \, |\, \max {\mathfrak e} =h)$,  the excursion $e$ can be split into two independent three-dimensional Bessel processes run until they hit $h$.   For $a\ge h$, and a three-dimensional Bessel process $R$ starting from $0$ stopped when hitting $h$, the It\^{o} formula together with the Dubins-Schwarz representation yield  that $\Phi(a-R)-\theta(a)$ is still a three dimensional  Bessel process run until it hits $\theta(a-h)-\theta(a)$ (this can also be seen as a special case of Theorem 2.2 equation (2.5) of Doney, Warren, Yor  \cite{dwy98} by taking $\alpha=0$ there).  It follows that under ${\mathfrak n}^+ (\cdot \, |\, \max {\mathfrak e} =h)$, $\Phi(a-{\mathfrak e})-\theta(a)$ is distributed as  ${\mathfrak e}$ under  ${\mathfrak n}^+ (\cdot \, |\, \max{\mathfrak e} =\theta(a-h)-\theta(a))$. Consequently, for any $a>0$,  \begin{eqnarray*}
 && 
 \int   {\mathfrak n}^+ (\d {\mathfrak e}) (1- e^{-f(\theta(a), \Phi(a-{\mathfrak e})-\theta(a) )}) 1_{\{ \max {\mathfrak e} < a \}}
 \\
 &=&
 \int_0^a  \frac{\d h}{2 h^2}  \int   (1- e^{-f(\theta(a), {\mathfrak e})}) {\mathfrak n}^+ \big(\d  {\mathfrak e} \, |\, \max {\mathfrak e} = \theta(a-h) - \theta(a)\big)
\\
&=&
 \frac{m^2}{(m+a)^2} \int_0^{|\theta(a)|} \frac{\d s}{2 s^2} \int   (1- e^{-f(\theta(a), {\mathfrak e})})  {\mathfrak n}^+ (\d  {\mathfrak e} \, |\, \max {\mathfrak e} = s)
\\
&=&
  \frac{m^2}{(m+a)^2}  \int {\mathfrak n}^+ (\d {\mathfrak e}) (1- e^{-f(\theta(a),   {\mathfrak e})}) 1_{\{ \max {\mathfrak e} < |\theta(a)|\}} ,
\end{eqnarray*}

\noindent  where the second equality follows from  a change of variables $s=\theta(a-h) - \theta(a)$.  It follows that \begin{eqnarray*}
&&
\e \Big[ e^{-\sum_{  \gamma\in{\mathcal L}_{\delta\over 2},  \min \gamma>0} f(  \min \Phi(\gamma)  ,   \Phi(\gamma)^\uparrow)}\Big] 
\\
&=&
\exp\Big(- \delta \int_0^\infty \d a   \frac{m^2}{(m+a)^2}  \int  {\mathfrak n}^+ (\d  {\mathfrak e}) (1- e^{-f(\theta(a),  {\mathfrak e})}) 1_{\{ \max {\mathfrak e} < |\theta(a)|\}}\Big)
\\
&=&
\exp\Big(- \delta  \int_0^m   \d y \int {\mathfrak n}^+ (\d  {\mathfrak e}) (1- e^{-f(-y, {\mathfrak e})}) 1_{\{ \max {\mathfrak e} < y\}} \Big),
\end{eqnarray*}

\noindent after a change of variables $y=|\theta(a)|$. This proves (a).

It remains to show (b).  Let $({\mathfrak e}_s, s>0)$ be the standard Brownian excursion process. It is well known that $(B_s, 0\le s \le \tau_\ell^B)$ can be constructed from $({\mathfrak e}_s, s\le \ell)$ (see Revuz and Yor \cite{revuz-yor} Chapter XII, Proposition 2.5).  Observe that  the process $\Phi(B_s, 0\le s \le \tau_\ell^B)$ can be constructed from $(\Phi({\mathfrak e}_s), s\le \ell)$ in the same way. To prove (b), it is enough to show that  $(\Phi({\mathfrak e}_s), s \le \ell)$ under the It\^{o} measure ${\mathfrak n}$, is distributed as $({\mathfrak e}_s, s \le \ell)$ under ${\mathfrak n}(\cdot \,|\, \inf_{s\le \ell} \min {\mathfrak e}_s >-m)$.  To this end, we use the same observation as in the proof of (a): for any $h>0$,  under ${\mathfrak n}^+ (\cdot \, |\, \max {\mathfrak e} =h)$, $\Phi({\mathfrak e})$ is distributed as  $-{\mathfrak e}$ under  ${\mathfrak n}^+ (\cdot \, |\, \max {\mathfrak e} =|\theta(h)|)$.  Consequently, for any  nonnegative measurable function $f$  on ${\mathcal K}$,  \begin{eqnarray*}
 \int  {\mathfrak n}^+ (\d {\mathfrak e}) (1- e^{-f(  \Phi({\mathfrak e}))})  
 &=&
 \int_0^\infty  \frac{\d h}{2 h^2}  \int   (1- e^{-f(-{\mathfrak e})}) {\mathfrak n}^+ (\d  {\mathfrak e} \, |\, \max {\mathfrak e} = |\theta(h)|)
\\
&=&
\int_0^m \frac{\d s}{2 s^2} \int   (1- e^{-f(-{\mathfrak e})})  {\mathfrak n}^+ (\d  {\mathfrak e} \, |\, \max {\mathfrak e} = s)
\\
&=&
 \int   {\mathfrak n}^+ (\d {\mathfrak e}) (1- e^{-f( -{\mathfrak e})}) 1_{\{ \max {\mathfrak e} <m\}} ,
\end{eqnarray*}

\noindent  where the second equality follows from  a change of variables $s=|\theta(h)|$.  It follows that \begin{eqnarray*}
 \int  {\mathfrak n} (\d {\mathfrak e}) (1- e^{-f(  \Phi({\mathfrak e}))})  
 &=& 
  \int  {\mathfrak n}^- (\d {\mathfrak e}) (1- e^{-f(-{\mathfrak e})})  + 
 \int  {\mathfrak n}^+ (\d {\mathfrak e}) (1- e^{-f( -{\mathfrak e})}) 1_{\{ \max {\mathfrak e} <m\}} 
 \\
 &=&
 \int  {\mathfrak n}(\d {\mathfrak e}) (1- e^{-f({\mathfrak e})}) 1_{\{ \min {\mathfrak e} >-m\}}  ,
\end{eqnarray*}

\noindent   which together with  the exponential formula for the excursion process,   yield that $(\Phi({\mathfrak e}_s), s \le \ell)$ under ${\mathfrak n}$  is distributed as $({\mathfrak e}_s, s \le \ell)$ under ${\mathfrak n}(\cdot \,|\, \inf_{s\le \ell} \min {\mathfrak e}_s >-m)$. This completes the proof of Theorem \ref{p:scaling1'}.  $\Box$

 \bigskip

\noindent {\large\bf Acknowledgements}

\medskip

We are grateful to Titus Lupu for stimulating discussions on the link between the PRBM and the Brownian loop  soup. We also thank an anonymous referee for useful suggestions on the paper. The project was partly supported by ANR MALIN.

\end{document}